\documentclass{article}

\usepackage[margin=1in]{geometry}
\usepackage[utf8]{inputenc}
\usepackage{amsmath, amssymb, tikz-cd, centernot, mathtools, csquotes,
            mathrsfs, hyperref, amsthm, xparse, microtype, bbm}
\usepackage[backend=biber, style=numeric, backref=true]{biblatex}
\usepackage[nameinlink,capitalize]{cleveref}
\usepackage{IEEEtrantools}

\hypersetup{
  colorlinks,
  linktoc=section,
  linkcolor=blue,
  citecolor=blue,
  urlcolor=blue,
}

\newcounter{smarttheorem:current} % the theorem most recently ended
\newcounter{smarttheorem:next} % the number to use for the next theorem we see

\NewDocumentCommand{\newsmarttheorem}{m o m o}{
  \IfNoValueTF{#2}{\newcounter{#1}}{}
  \newtheorem{hidden#1}[\IfValueTF{#2}{#2}{#1}]{#3}

  \NewDocumentEnvironment{#1}{o}{
    \IfNoValueTF{##1}{\begin{hidden#1}}{\begin{hidden#1}[{##1}]}
      \label{smarttheorem:\arabic{smarttheorem:next}}
      \edef\currentlabel{\arabic{smarttheorem:next}}
      \stepcounter{smarttheorem:next}
  }{
      \setcounter{smarttheorem:current}{\currentlabel}
    \end{hidden#1}
  }
  \Crefname{hidden#1}{#3}{\IfNoValueTF{#4}{#3s}{#4}}
}

\let \vqedsymbol \qedsymbol

\NewDocumentEnvironment{pf}{o o}{
  \IfValueTF{#2}{
    \edef \qedsymbollabel {#2}
    \renewcommand \qedsymbol {\vqedsymbol \, \qedsymbollabel}
  }{
    \renewcommand \qedsymbol \vqedsymbol
  }
  \IfValueTF{#1}{\begin{proof}[#1]}{\begin{proof}}
}{
  \end{proof}
}

\NewDocumentEnvironment{rpf}{O{\proofname}
O{smarttheorem:\arabic{smarttheorem:current}}}{
  \begin{pf}[#1][\noexpand\Cref{#2}]
}{
  \end{pf}
}

\NewDocumentEnvironment{lrpf}{m}{
  \begin{rpf}[Proof of \cref{#1}][#1]
}{
  \end{rpf}
}

\newcommand{\pss}{\item[\((\subseteq)\)]}
\newcommand{\psps}{\item[\((\supseteq)\)]}

\NewDocumentEnvironment{ea}{O{rCl}}{
  \begin{IEEEeqnarray*}{#1}
}{
  \end{IEEEeqnarray*}
  \ignorespacesafterend
}

\NewDocumentEnvironment{tcd}{s}{
  \IfBooleanTF{#1}{\begin{equation}}{\begin{equation*}}
    \begin{tikzcd}
}{
    \end{tikzcd}
  \IfBooleanTF{#1}{\end{equation}}{\end{equation*}}
  \ignorespacesafterend
}

\newcounter{case}
\newenvironment{caselist}{
  \setcounter{case}{0}\begin{description}
}{
  \end{description}
}

\newcommand \case {\setcounter{case}{\value{case}+1}\item[Case \thecase.]}

\numberwithin{theoremcounter}{section}
\theoremstyle{plain}
\newsmarttheorem{theorem}[theoremcounter]{Theorem}
\newsmarttheorem{lemma}[theoremcounter]{Lemma}
\newsmarttheorem{proposition}[theoremcounter]{Proposition}
\newsmarttheorem{claim}[theoremcounter]{Claim}
\newsmarttheorem{subclaim}[theoremcounter]{Subclaim}
\newsmarttheorem{corollary}[theoremcounter]{Corollary}[Corollaries]
\newsmarttheorem{fact}[theoremcounter]{Fact}
\newsmarttheorem{todo}{TODO}

\newtheorem*{theorem-non-num}{Theorem}
\newtheorem*{corollary-non-num}{Corollary}

\theoremstyle{definition}
\newsmarttheorem{definition}[theoremcounter]{Definition}
\newsmarttheorem{notation}[theoremcounter]{Notation}

\theoremstyle{remark}
\newsmarttheorem{question}[theoremcounter]{Question}
\newsmarttheorem{remark}[theoremcounter]{Remark}
\newsmarttheorem{scratch}[theoremcounter]{Rough work}
\newsmarttheorem{aside}[theoremcounter]{Aside}
\newsmarttheorem{conjecture}[theoremcounter]{Conjecture}
\newsmarttheorem{exercise}[theoremcounter]{Exercise}
\newsmarttheorem{example}[theoremcounter]{Example}
\newsmarttheorem{ednote}[theoremcounter]{Editor's note}

\renewcommand \epsilon \varepsilon
\renewcommand \phi \varphi

\renewcommand \th {^\mathrm{th}}
\newcommand \telse {\text{else}}
\newcommand \tif {\text{if }}

\newcommand \pip [2] {\item[\underline{(#1)$\implies$(#2)}]}

\newcommand \symdif {\mathbin\triangle}

\def\Indep#1#2{#1\setbox0=\hbox{$#1x$}\kern\wd0\hbox to 0pt{\hss$#1\mid$\hss}
\lower.9\ht0\hbox to 0pt{\hss$#1\smile$\hss}\kern\wd0}

\def\notindep#1#2{#1\setbox0=\hbox{$#1x$}\kern\wd0
\hbox to 0pt{\mathchardef\nn=12854\hss$#1\nn$\kern1.4\wd0\hss}
\hbox to
0pt{\hss$#1\mid$\hss}\lower.9\ht0 \hbox to 0pt{\hss$#1\smile$\hss}\kern\wd0}

\newcommand \tfae {the following are equivalent}
\newcommand \Tfae {The following are equivalent}

\newenvironment{menumerate}{\mbox{}\begin{enumerate}}{\end{enumerate}}

\DeclarePairedDelimiter \abs \lvert \rvert

\DeclarePairedDelimiter \pars ( )
\DeclarePairedDelimiter \bracks [ ]

\DeclarePairedDelimiter \floor \lfloor \rfloor

\DeclarePairedDelimiterX \set [1] \{ \} {\, #1 \,}

\newcommand \col [2] {\begin{pmatrix}#1\\#2\end{pmatrix}}
\newcommand \tcol [3] {\begin{pmatrix}#1\\#2\\#3\end{pmatrix}}
\newcommand \fcol [4] {\begin{pmatrix}#1\\#2\\#3\\#4\end{pmatrix}}
\newcommand \powmult [1] {\times \mathord\restriction #1^\mathbb{N}}
\newcommand \dpowmult {\powmult{d}}

\DeclareMathOperator \id {id}

\DeclareMathOperator \thy {Th}

\title{Automata and tame expansions of \((\mathbb{Z},+)\)}
\author{Christopher Hawthorne\thanks{%
This work was partially supported by an NSERC PGS-D and an NSERC CGS-D.}}

\begin{document}
\maketitle
\begin{abstract}
  The problem of characterizing which automatic sets of integers are stable is
  here solved. Given a positive integer \(d\) and a subset \(A\subseteq
  \mathbb{Z}\) whose set of representations base \(d\) is recognized by a
  finite automaton, a necessary condition is found for \(x+y\in A\) to be a
  stable formula in \(\thy(\mathbb{Z},+,A)\).  Combined with a theorem of Moosa
  and Scanlon this gives a combinatorial characterization of the
  \(d\)-automatic \(A\subseteq \mathbb{Z}\) such that \((\mathbb{Z},+,A)\) is
  stable. This characterization is in terms of what were called
  \emph{\(F\)-sets} in \cite{moosa04} and \emph{elementary \(p\)-nested sets}
  in \cite{derksen07}. Automata-theoretic methods are also used to produce some
  NIP expansions of \((\mathbb{Z},+)\), in particular the expansion by the
  monoid \((d^\mathbb{N},\times )\).
\end{abstract}
\section{Introduction}
In \cite{palacin18} Palacín and Sklinos give examples of stable expansions of
\(\thy(\mathbb{Z},+)\), and pose the following general question:
\begin{question}
  \label{question:stable}
  For which \(A\subseteq \mathbb{Z}\) is \(\thy(\mathbb{Z},+,A)\) stable?
\end{question}
The project of finding sufficient topological or algebraic conditions on \(A\)
for stability has been taken up in other recent work; see for example
\cite{conant19} and \cite{lambotte20}. The theme also appeared some fifteen
years earlier: the results of Moosa and Scanlon in \cite{moosa04} imply that
\((\mathbb{Z},+,A)\) is stable whenever \(A\) is an \emph{\(F\)-set} (see
\cref{def:fsets}). This includes for example the case \(A = d^\mathbb{N}\),
whose stability was rediscovered in \cite{palacin18}.

In this paper, we consider \cref{question:stable} when \(A\subseteq
\mathbb{Z}\) is a \(d\)-automatic set for some \(d\ge 2\). Automatic sets are
reviewed in \cref{sec:prelim}, but let us recall here informally that this
means there is a finite machine that takes strings of digits from
\(\set{-d+1,\ldots ,d-1}\) as input and accepts exactly those strings that are
representations base \(d\) of an element of \(A\).

Instead of asking when the first-order theory of \((\mathbb{Z},+,A)\) is
stable, we will focus on a local, and hence combinatorial, notion of stability
which we now briefly recall. If \(R\subseteq X\times X\) is a binary relation
on a set \(X\) then an \emph{\(N\)-ladder} for \(R\) is some \(a_0,\ldots
,a_N,b_0,\ldots ,b_N\in X\) such that \((a_i,b_j)\in R\) if and only if \(i\le
j\). The relation \(R\) is \emph{\(N\)-stable} if there is no \(N\)-ladder for
\(R\), and is \emph{stable} if it is \(N\)-stable for some \(N\). If \((G,+)\)
is a group and \(A\subseteq G\) we say that \(A\) is \emph{stable in \(G\)} if
\(x+y\in A\) is a stable binary relation on \(G\).\footnote{It is worth noting
  that this terminology conflicts somewhat with the terminology used by Conant
  in \cite{conant19}, in which he calls \(A\subseteq \mathbb{N}\) ``stable'' if
  \(\thy(\mathbb{Z},+,A)\) is stable. His is a stronger notion than ours, which
  is equivalent to saying that \(\phi (x;y)\) given by \(x+y\in A\) is a stable
  formula in \(\thy(G,+,A)\).
}

Here is our main result, which appears as \cref{thm:charstableauto} below.
\begin{theorem-non-num}
  Suppose \(A\) is \(d\)-automatic and stable in \((\mathbb{Z},+)\). Then \(A\)
  is a finite Boolean combination of
  \begin{itemize}
    \item
      cosets of subgroups of \((\mathbb{Z},+)\), and
    \item
      translates of finite sums of sets of the form
      \[C(a;\delta ):= \set{a +
        d^\delta a+\cdots + d^{n\delta }a : n < \omega }
      \]
      where \(a\in\mathbb{Z}\) and \(\delta \) is a positive integer.
  \end{itemize}
\end{theorem-non-num}
These sets are of Diophantine significance in positive characteristic,
appearing in both the isotrivial Mordell-Lang theorem of \cite{moosa04} and the
Skolem-Mahler-Lech theorem of \cite{derksen07}; see \cite[Section 3]{bell19}
for an account of the connection with the latter.  In particular, combining our
main theorem with the results of \cite{moosa04} yields:
\begin{corollary-non-num}
  Suppose \(A\subseteq \mathbb{Z}\) is \(d\)-automatic. Then \tfae{}:
  \begin{enumerate}
    \item
      \(\thy(\mathbb{Z},+,A)\) is stable.
    \item
      \(A\) is stable in \((\mathbb{Z},+)\).
    \item
      \(A\) is a finite Boolean combination of cosets of subgroups of
      \((\mathbb{Z},+)\) and translates of finite sums of sets of the form
      \(C(a;\delta )\).
  \end{enumerate}
\end{corollary-non-num}
This appears as \cref{cor:autosparsechar} below.

Automatic sets separate naturally into sparse and non-sparse sets, with
``sparse'' meaning that the number of accepted strings grows polynomially with
length---see \cref{def:dsparse} for a precise formulation.  The first case of
the main theorem that we consider is when \(A\) is \(d\)-sparse. In fact, here
we can work more generally in Cartesian powers of \((\mathbb{Z},+)\). So, in
\cref{thm:charstablesparse} below we prove that if \(A\subseteq \mathbb{Z}^m\)
is \(d\)-sparse and stable in \((\mathbb{Z}^m,+)\) then it is a finite Boolean
combination of translates of finite sums of sets of the form \(C(a;\delta )\)
where now \(a\) is an element of \(\mathbb{Z}^m\).

We then turn our attention to \(d\)-automatic sets that aren't \(d\)-sparse.
We show that for \(A\subseteq \mathbb{Z}\) \(d\)-automatic
but not \(d\)-sparse, if \(A\) is stable in \((\mathbb{Z},+)\) then \(A\) is
generic (i.e.\ finitely many translates cover \(\mathbb{Z}\)). This is
\cref{thm:stablenongen}. In particular, every \(d\)-automatic subset of
\(\mathbb{N}\) that is stable in \((\mathbb{Z},+)\) must be \(d\)-sparse.
Actually, this consequence of our \cref{thm:stablenongen} can also be deduced
by combining \cite[Theorem 8.8]{conant19} and \cite[Theorem 1.1]{bell18}, but
the general statement requires significantly more work.

\Cref{thm:charstablesparse,thm:stablenongen}, together with some stable group
theory, yield the main theorem.

In a somewhat different direction, we conclude the paper by using
automata-theoretic methods to produce two NIP expansions of \((\mathbb{Z},+)\):
namely \((\mathbb{Z},+,<,d^\mathbb{N})\) and
\((\mathbb{Z},+,d^\mathbb{N},\dpowmult)\), in
\cref{thm:presbpowersnip,thm:multpowersnip}, respectively. The former was
recently obtained by Lambotte and Point in \cite{lambotte20} using different
methods, but the latter is a new example.

\subsection*{Acknowledgements}
I am very grateful to Gabriel Conant, who, upon viewing an earlier draft of
this paper, pointed out to me that the cases dealt with in
\cref{thm:charstablesparse,thm:stablenongen} were sufficient to prove the
general case. I am also grateful to Jason Bell, in conversations with whom the
main theorem was first articulated as a conjecture. Finally, I am deeply
grateful to my advisor, Rahim Moosa, for excellent guidance, thorough editing,
and many helpful discussions.

\section{Preliminaries on automatic sets}
\label{sec:prelim}
We briefly recall regular languages and finite automata; see \cite{yu97} for a
more detailed presentation.
\begin{definition}
  For a finite set \(\Lambda \) viewed as an \emph{alphabet} we denote by
  \(\Lambda ^*\) the set of \emph{words} over \(\Lambda \), namely finite
  strings of letters from \(\Lambda \).  The class of \emph{regular languages}
  over \(\Lambda \) is the smallest subset of \(\mathcal{P}(\Lambda ^*)\) that
  contains all finite sets and is such that if \(A,B\) are regular then so are
  \(A\cup B\), \(AB\), and \(A^*\).
\end{definition}
A \emph{deterministic finite automaton} (DFA) over a finite alphabet \(\Lambda
\) is a tuple \(\mathcal{A} = (Q,q_0,F,\delta )\) where \(Q\) is a finite set
of states, \(q_0\in Q\) is the start state, \(F\subseteq Q\) is the set of
finish states, and \(\delta \colon Q\times \Lambda \to Q\) is the transition
function: if \(\mathcal{A}\) is in state \(q\in Q\) and reads the letter
\(\ell \in\Lambda \) then it moves into state \(\delta (q,\ell )\).  We
identify \(\delta \) with its natural extension \(Q\times \Lambda ^*\to Q\)
inductively by \(\delta (q,\ell _1\cdots \ell _{n+1}) = \delta (\delta
(q,\ell _1\cdots \ell _n),\ell _{n+1})\).  The set \emph{recognized}
by \(\mathcal{A}\) is \(\set{\sigma \in\Lambda ^* : \delta (q_0, \sigma)\in
F}\). A fundamental fact (see \cite[Lemma 2.2, Section 3.2, and Section
3.3]{yu97}) is that the regular languages are precisely the sets recognized by
DFAs.

It turns out some behaviour of automata can be captured in Presburger
arithmetic:
\begin{proposition}
  \label{prop:sparseintersregpresb}
  Suppose \(\Lambda \) is an alphabet; suppose \(L\subseteq \Lambda ^*\) is
  regular and \(\sigma _1,\ldots ,\sigma _n\in\Lambda
  ^*\). Then \(\set{(t_1,\ldots ,t_n)\in\mathbb{N}^n : \sigma _1^{t_1}\cdots
  \sigma _n^{t_n}\in L}\) is definable in \((\mathbb{N},+)\).
\end{proposition}
\begin{rpf}
  Fix an automaton \((Q,q_0,F,\delta )\) for \(L\).  We apply induction on
  \(n\) to show that \(\delta (q_1,\sigma _1^{t_1}\cdots \sigma _n^{t_n}) =
  q_2\) is definable in \((\mathbb{N},+)\) for all \(q_1,q_2\in Q\). The case
  \(n=0\) is vacuous.  For the induction step, note since there are finitely
  many states that \(\delta (q_1, \sigma _1^t ) \) is eventually cyclic in \(t
  \); so there is \(N\) such that for \(t \ge N\) we have that \(\delta
  (q_1,\sigma _1^t ) \) depends only on the congruence class of \(t \)
  modulo some \(\mu \). But then our set is defined by
  \[\bigvee _{t < N}\pars[\Big]{(t _1=t )\wedge \delta (\delta
    (q_1,\sigma _1^t ), \sigma _2^{t _2}\cdots \sigma _n^{t _n}) =
    q_2}
    \vee 
    \bigvee _{ t  < \mu } \pars[\Big]{(t _1\in N+t+\mu
      \mathbb{N}) \wedge \delta (\delta (q_1,\sigma _1^{N+t} ), \sigma _2^{t
    _2}\cdots \sigma _n^{t _n})= q_2}
  \]
  which is definable in \((\mathbb{N},+)\) by the induction hypothesis.

  But then \(\set{(t_1,\ldots ,t_n) \in\mathbb{N}^n : \sigma _1^{t_1}\cdots
  \sigma _n^{t_n}\in L}\) is the union over \(q\in F\) of \(\set{(t_1,\ldots
    ,t_n)\in\mathbb{N}^n : \delta (q_0, \sigma _1^{t_1}\cdots \sigma _n^{t_n})
  = q}\), which is definable in \((\mathbb{N},+)\).
\end{rpf}
We are primarily interested in the case where the strings in question are
representations of integers. Fix a positive integer \(d\). Evaluating a
string base \(d\) gives a natural map \([\cdot ]\colon
\mathbb{Z}^*\to\mathbb{Z}\) via
\[[k_0k_1\cdots k_n] = \sum_{i=0}^n k_id^i.
\]
Note that unlike usual base \(d\) representations the most significant digit
occurs last, not first.
\begin{definition}
  We let \(\Sigma =\set{0,\ldots ,d-1}\) and \(\Sigma _\pm = \set{-d+1,\ldots
  ,d-1}\).  We say \(A\subseteq \mathbb{Z}\) is \emph{\(d\)-automatic} if
  \(\set{\sigma \in\Sigma _\pm^* : [\sigma ]\in A}\) is a regular language over
  \(\Sigma _\pm\).
\end{definition}
There is a natural extension of this notion to \(\mathbb{Z}^m\) for \(m\ge 1\).
We first extend our map \([\cdot ]\) to \((\mathbb{Z}^m)^*\to\mathbb{Z}^m\): we
set
\[\bracks*{\tcol{k_{10}}\vdots {k_{m0}}\cdots \tcol{k_{1n}}\vdots {k_{mn}}}
  = \tcol{\relax[k_{10}\cdots k_{m0}]}\vdots {\relax [k_{1n}\cdots k_{mn}]}
\]
\begin{definition}
  We say \(A\subseteq \mathbb{Z}^m\) is \emph{\(d\)-automatic} if
  \(\set{\sigma \in(\Sigma _\pm ^m)^*: [\sigma ]\in A}\) is a regular
  language over \(\Sigma _\pm^m\).
\end{definition}
A note on exponential notation: we use \(\Lambda ^m\) to denote the alphabet
\(\Lambda \times \cdots \times \Lambda \). This contrasts its usual meaning in
formal languages, namely the set of words over \(\Lambda \) of length \(m\); we
use \(\Lambda ^{(m)}\) to denote this. We will use \(\sigma ^n\) to denote the
\(n\)-fold concatenation of \(\sigma \) with itself; it should be clear from
context whether an instance of exponential notation refers to iterated string
concatenation or iterated multiplication.

Of course different strings can represent the same integer. It is useful to fix
a canonical representation.
\begin{definition}
  The \emph{canonical representation} of \(0\) is the empty word \(\epsilon \).
  The \emph{canonical representation} of a positive integer \(a\) is its usual
  representation base \(d\) in \(\Sigma ^*\) (though with the order reversed).
  The \emph{canonical representation} of a negative integer \(a\) is
  \((-k_0)\cdots (-k_n)\) where \(k_0\cdots k_n\) is the canonical
  representation of \(-a\). Finally, the \emph{canonical representation} of a
  tuple \(\tcol{a_1}\vdots {a_m}\) is \(\tcol{k_{10}}\vdots {k_{m0}}\cdots
  \tcol{k_{1n}}\vdots {k_{mn}} \) where \(n+1\) is the maximum of the lengths
  of the canonical representations of the \(a_i\), and \(k_{i0}\cdots k_{in}\)
  is the canonical representation of \(a_i\) for each \(i\), possibly padded
  with trailing zeroes to make them of length \(n+1\). 
\end{definition}
Note that the canonical representation of an integer base \(d\) is a word over
\(\Sigma _\pm\), and if the integer happens to be non-negative then it is a
word over \(\Sigma \).
\begin{example}
  The canonical representation base \(10\) of \(\col{-23}{432}\) is
  \(\col{-3}{2}\col{-2}{3}\col{0}{4}\).
\end{example}
\begin{remark}
  \label{rem:autorobust}
  Automaticity is robust under changes in the allowed representations.  Indeed,
  from \cite[Proposition 7.1.4]{frougny02} (together with some basic tools of
  automata theory---see e.g.\ \cite[Theorems 2.16 and 4.2]{yu97}) one can show
  that the following are both equivalent to \(A\) being \(d\)-automatic:
  \begin{enumerate}
    \item
      The set of canonical representations of elements of \(A\) is regular.
      Note that this is essentially the definition given in \cite[Section
      5.3]{adamczewski12}. In particular, our definition generalizes the
      classical notion of \(d\)-automatic subsets of \(\mathbb{N}\) (see e.g.\
      \cite[Chapter 5]{allouche03}).
      % (See e.g.\ \cite[Chapter 5]{allouche03}, though note that they use
      % most-significant-digit-first representations of integers, which is the
      % reverse of our convention. This doesn't affect the definition of
      % \(d\)-automaticity, though, since \(L\) is regular if and only if
      % \(\set{\sigma ^R : \sigma \in L}\) is, where \(\sigma ^R\) is the
      % reversal of \(\sigma \); see \cite[Corollary 4.3.5]{allouche03}.)
    \item
      For some (equivalently all) finite \(\Lambda \subseteq \mathbb{Z}^m\)
      containing \(\Sigma _\pm^m\), 
      \(\set{\sigma \in\Lambda ^* : [\sigma ]\in A}\) is regular over
      \(\Lambda\).
  \end{enumerate}
\end{remark}
We are particularly interested in sparsity among \(d\)-automatic sets.  If
\(\Lambda \) is an alphabet we say \(L\subseteq \Lambda ^*\) is \emph{sparse}
if it is regular and the map \(k\mapsto \abs{\set{\sigma \in L : \abs\sigma \le
k}}\) grows polynomially in \(k\). Several equivalent formulations of sparsity
are known; we will in particular make use of the following characterization:
\begin{fact}[\cite[Proposition 7.1]{bell19}]
  \label{fact:sparsechar}
  If \(L\subseteq \Lambda ^*\) then \(L\) is sparse if and only if it is a
  finite union of sets of the form \(u_0w_1^*u_1\cdots u_{n-1}w_n^*u_n =
  \set{u_0w_1^{r_1}u_1\cdots w_n^{r_n}u_n : r_1,\ldots ,r_n\ge 0}\) for some
  \(u_0,\ldots ,u_n,w_1,\ldots ,w_n\in\Lambda ^*\).
\end{fact}
\begin{definition}
  \label{def:dsparse}
  We say \(A\subseteq \mathbb{Z}^m\) is \emph{\(d\)-sparse} if the set of
  canonical representations base \(d\) of elements of \(A\) is a sparse
  language over \(\Sigma _\pm^m\).
\end{definition}
Note by \cref{rem:autorobust} that \(d\)-sparse sets are \(d\)-automatic. In
fact, \(d\)-sparsity is equivalent to the existence of some finite \(\Lambda
\supseteq \Sigma _\pm^m\) and some sparse \(L\subseteq \Lambda ^*\) such that
\(A = [L]\), but we will not need this.
\section{Characterizing stable sparse sets}
\label{sec:sparsestable}
In this section we give a complete characterization of the \(d\)-sparse subsets
of \(\mathbb{Z}^m\) that are stable in \((\mathbb{Z}^m,+)\). Observe that not
every \(d\)-sparse set is stable: assuming \(d>2\) the set \(A = [0^*10^*2]\)
is \(d\)-sparse by \cref{fact:sparsechar}, but it is not hard to verify that
\(d^i + 2d^{j+1}\in A\) if and only if \(i \le j\). 

The main result of this section:
\begin{theorem}
  \label{thm:charstablesparse}
  Suppose \(A\subseteq \mathbb{Z}^m\) is \(d\)-sparse. If \(A\) is stable in
  \((\mathbb{Z}^m,+)\) then \(A\) is a finite Boolean combination of translates
  of finite sums of sets of the form
  \[C(a ; \delta ) = \set{a + d^\delta a + \cdots + d^{n\delta }a : n < \omega
    }
  \]
  where \(a\in\mathbb{Z}^m\) and \(\delta >0\).
\end{theorem}
Sets of this form, namely finite unions of translates of finite sums of
\(C(a;\delta )\), were studied in \cite{moosa04} in a more general setting,
where they were called ``groupless \(F\)-sets''; here \(F\) is the
multiplication-by-\(d\) endomorphism on \((\mathbb{Z}^m,+)\). When \(m=1\),
these groupless \(F\)-sets were rediscovered in a different context by Derksen
\cite{derksen07} as ``elementary \(p\)-nested'' subsets of \(\mathbb{Z}\); see
\cite[Proposition 3.2]{bell19} for a proof that they agree up to finite
symmetric differences.

Combined with the results of \cite{moosa04} this theorem yields complete
answers to both \cref{question:stable,question:localstable} for \(d\)-sparse
subsets of \(\mathbb{Z}^m\). See \cref{cor:sparsestablechar} below.

Before proving the theorem let us make some observations that may give the
reader a better feel for the automata-theoretic nature of the sets \(C(a;\delta
)\).
\begin{remark}
  \label{rem:fsets}
  \begin{menumerate}
    \item
      \(d^\mathbb{N}=\set1\cup(1+C(d-1;1))\).
    \item
      If \(a = [\sigma ]\) where \(\sigma \in (\mathbb{Z}^m)^*\) is of
      length \(\delta \) then \(C(a;\delta ) = \set{[\sigma ^n]:n>0}\).
    \item
      Every translate of a finite sum of sets of the form \(C(a;\delta )\) is
      \(d\)-sparse.
    \item
      Let \(\mathcal{C}\) denote the collection of subsets of \(\mathbb{Z}\) of
      the form \(b+C(a;\delta )\) for some \(a,b\in\mathbb{Z}\) and \(\delta
      >0\). Let \(\mathcal{E}\) be the collection of subsets of \(\mathbb{Z}\)
      of the form \([uv^*w]\) for \(u,v,w\in \Sigma ^*\) or
      \(u,v,w\in(-\Sigma)^*\).  Then up to finite symmetric differences,
      \(\mathcal{C}\) and \(\mathcal{E}\) agree.
  \end{menumerate}
\end{remark}
\begin{rpf}
  Parts (1) and (2) are easily verified by hand. Part (3) is observed in
  \cite{bell19}; see the proof of Theorem 7.4 therein.  We prove (4).
  \begin{description}
    \psps
      We will see in the proof of \cref{lemma:sparseincycles} below that we can
      write \([uv^*w]\) as a translate of \([\tau ^*]\) for some \(\tau \in
      \mathbb{Z}^*\).  But then by part (2) this has finite symmetric
      difference with \(C([\tau ]; \abs\tau )\).
    \pss
      Suppose first that we are given \(C(a;\delta )\); by negating if
      necessary we may assume \(a\ge 0\).
      Pick a representation \(\sigma \in \mathbb{Z}^*\) of \(a\) of length
      \(\delta \); then by part (2) we are interested in the canonical
      representations of \([ \sigma ^*\setminus \set\epsilon ]\).  For
      \(0<i<\omega \) write \([\sigma ^i] = b_i d^{i\delta } + c_i\) where
      \(0\le c_i < d^{i\delta }\); so \(b_i\) is the ``carry'' when adding up
      \(a + F^\delta a +\cdots +F^{(i-1)\delta }a\) and cutting off after
      \(i\delta \) digits.  Then
      \[[\sigma ^{i+1}] = d^\delta [\sigma ^i] + [\sigma ] = b_id^{(i+1)\delta
        } + d^\delta c_i + [\sigma ]\ge b_id^{(i+1)\delta }
      \]
      so \(b_{i+1}\ge b_i\). But \([\sigma ^i] = a\frac{d^{i\delta
      }-1}{d^\delta -1} \le ad^{i\delta }\); so each \(b_i\le a\). So the
      \(b_i\) are eventually constant, say \(b_N = b_{N+1} = \cdots \). Let \(p
      = b_N+[\sigma ]\bmod d^\delta \). Then
      \[[\sigma ^{N+k+1} ]
        = [\sigma ^{N+k}] + d^{(N+k)\delta }[\sigma ]
        = b_{N+k}d^{(N+k)\delta } + c_{N+k} + d^{(N+k)\delta }[\sigma ]
        = d^{(N+k)\delta } (b_N + [\sigma ]) + c_{N+k}
      \]
      so
      \[c_{N+k+1} = [\sigma ^{N+k+1}]\bmod d^{(N+k+1)\delta }
        = d^{(N+k)\delta } p + c_{N+k}
      \]
      (since \(c_{N+k} < d^{(N+k)\delta }\)).  So inductively we get \(c_{N+k}
      = c_N + d^{N\delta }p + d^{(N+1)\delta }p + \cdots + d^{N+k-1}p\). So
      \[[\sigma ^{N+k}] = b_N d^{(N+k)\delta } + c_N + d^{N\delta }p + \cdots +
        d^{N+k-1}p
      \]
      So if \(u\in\Sigma ^{(N\delta) },v\in \Sigma ^{(\delta)} , w\in\Sigma ^*\)
      represent \(c_N, p, b_N\) respectively then \([\sigma ^{N+k}] =
      [uv^kw]\). So \(C(a ;\delta ) = [\sigma ^*\setminus \set\epsilon ]\) has
      finite symmetric difference with \([uv^*w]\), as desired.

      It remains to show that a translate of a single cycle takes the desired
      form; by above it suffices to show that a translate of \([uv^*w]\), say
      by \(\gamma \in\mathbb{Z}\), has finite symmetric difference from some
      \([xy^*z]\).  (Again we may assume \([uv^*w]\subseteq \mathbb{N}\).) If
      \(u,v\in (d-1)^*\) then \(\gamma + [uv^*w] = (\gamma -1) + [0^{\abs
      u}(0^{\abs v})^*\tau ]\) where \([\tau ] = [w] + 1\); so we may assume
      \(uv\notin (d-1)^*\). So for some \(N\) we get that \(0\le \gamma +
      [uv^N] < d^{\abs{uv^N}}\); so if \(\sigma \in \Sigma ^{(\abs{uv^N})}\) has
      \([\sigma ] = \gamma +[uv^N]\) then \(\gamma +[uv^{N+k}w] = [\sigma
      v^kw]\). So \(\gamma +[uv^*w]\) has finite symmetric difference from
      \([\sigma v^*w]\).
      \qedhere
  \end{description}
\end{rpf}
We begin working towards a proof of \cref{thm:charstablesparse}. Our approach
requires that we first understand the stable formulas in
\((\mathbb{N},0,S,\delta \mathbb{N},<)\) where \(S\) is the successor function
and \(\delta \) is a fixed positive integer. The following proposition, which
is of independent interest, is likely known, but as we could find no reference
we include a proof here for completeness.
\begin{proposition}
  \label{prop:stabledeltaeqn}
  Fix \(\thy(\mathbb{N},0,S,\delta \mathbb{N},<)\) as the ambient theory. Let
  \(L_\delta =\set{0,S,P_\delta }\) and \(L_{\delta ,<} = L_\delta \cup\set<\).
  Suppose \(\phi (x_1,\ldots ,x_n)\in L_{\delta ,<}\) is quantifier-free and
  stable with respect to any partition of the variables. Then \(\phi \) is
  equivalent to a quantifier-free \(L_\delta \)-formula.\footnote{
    In fact both \(\thy(\mathbb{N},0,S,\delta \mathbb{N},<)\) and
    \(\thy(\mathbb{N},0,S,\delta \mathbb{N})\) admit quantifier elimination;
    however we don't make use of this in either the proof or the application of
    this proposition.
  }
\end{proposition}
\begin{rpf}
  We apply induction on \(n\); the case \(n=0\) is vacuous.

  With an eye towards constraining the atomic subformulas of \(\phi \), we
  rewrite \(\phi \) as follows:
  \begin{itemize}
    \item
      Replace any occurrence of \(S^ex_i < K\) by a disjunction of equalities
      in the obvious way, and of \(S^ex_i < S^f x_j\) for \(e \ge f\) by
      \(S^{e-f}x_i < x_j\). Using this and the fact that \(t_1\le t_2\iff \neg
      (t_1> t_2)\), we may assume all atomic inequalities take the form \(S^e
      x_i < x_j\).
    \item
      Replace any occurrence of \(S^ex_i = K\) by \(x_i = K-e\) and of \(S^ex_i
      = S^f x_j\) for \(e\ge f\) by \( S^{e-f}x_i=x_j\).
  \end{itemize}
  So we may assume the atomic subformulas of \(\phi \) take the following
  forms:
  \begin{itemize}
    \item
      \(x_i\equiv K\pmod \delta \)
    \item
      \(x_i = K\)
    \item
      \(S^ex_i < x_j\)
    \item
      \(S^ex_i=x_j\)
  \end{itemize}
  Let \(M\) be greater than both the largest \(K\) appearing in \(\phi \) and
  the largest \(e\) with \(S^e\) appearing in \(\phi \). Note that the truth
  value of \(\phi (\overline{a})\) is determined by the truth value of the
  above formulas on \(\overline{a}\). Furthermore since we may assume in said
  formulas that \(K<M\) and \(e<M\), we get that there are finitely many such
  formulas; call the set of such formulas \(\Delta \). So we can write \(\phi
  \) as a finite disjunction of consistent conjunctions of the form
  \[\psi _f = \bigwedge _{\theta \in\Delta }\theta ^{f(\theta )}
  \]
  for some \(f\colon \Delta \to\set{0,1}\). (Here \(\psi _f^0\) denotes
  \(\neg\psi _f\) and \(\psi _f^1\) denotes \(\psi _f\).)

  Fix one such disjunct \(\psi _f\); we will show that \(\psi _f\) implies
  some \(L_\delta  \)-formula that in turn implies \(\phi \).
  \begin{caselist}
    \case
      Suppose \(\psi _f\) contains a conjunct of the form
      \(S^ex_{j_1}=x_{j_2}\) or \(K = x_{j_2}\). Define a term \(t\) to be
      \(S^ex_{j_1}\) in the former case and \(K\) in the latter case, and
      consider \(\phi '(\overline{x'}) = \phi '(x_1,\ldots
      ,x_{j_2-1},x_{j_2+1},\ldots ,x_n)\) obtained by substituting \(x_{j_2} =
      t\) into \(\phi \); this is stable because \(\phi \) is, and because
      \(t\) involves at most one of the \(x_i\). It also contains one fewer
      variable, so by the induction hypothesis is equivalent to a
      quantifier-free \(L_\delta \)-formula \(\theta (\overline{x'}) \).  But
      then \(\theta (\overline{x'})\wedge (x_{j_2}=t)\) is our desired formula.
    \case
      Suppose \(\psi_f \) contains no such conjuncts; let \(\psi
      _0(\overline{x})\) be the conjunction of the negations of such.
      Examining \(\Delta \) we see that since \(\psi_f \) is consistent it must
      take the form
      \[\underbrace{\psi _0(\overline{x})\wedge \pars*{\bigwedge _{i=1}^n
        x_i\equiv K_i\pmod \delta }}_\chi \wedge ( x_{\sigma (1)}<\cdots
        <x_{\sigma (n)})
      \]
      for some \(K_i<\delta \) and some \(\sigma \in S_n\). (Note that formulas
      of the form \(S^ex_i<x_j\) for \(e<M\) are implied by \(x_i < x_j\) and
      \(\psi _0(\overline{x})\), so we may safely omit them.) I claim that
      \(\chi \) is our desired formula. It is clear that \(\models \psi_f
      \rightarrow \chi \); it remains to show that \(\models \chi \rightarrow
      \phi \). Fix \(1<j\le n\), and suppose for contradiction we had \(\models
      \neg\phi (\overline{a})\) for some
      realization \(\overline{a}\) of \(\chi \wedge (x_{((j-1\;j)\sigma)(1)
      }<\cdots <x_{((j-1\;j)\sigma )(n)})\); note this last formula takes the
      form \(\psi _{f'}\) for some \(f'\colon \Delta \to\set{0,1}\). So by
      definition of \(\Delta \) we get that \(\models \neg\phi
      (\overline{a})\) for all realizations \(\overline{a}\) of \(\psi _{f'}\).

      I claim this implies that \(\phi(x_1,\ldots ,x_{j-1}; x_j,\ldots ,x_n) \)
      has the order property. Indeed, fix \(N<\omega \); we construct a ladder
      of length \(N\). For clarity we assume \(\sigma = \id\); the argument
      generalizes with little effort.  Pick \(a_1\ge M\) such that \(a_1\equiv
      K_1\pmod\delta \), and inductively pick \(a_{i+1}\ge a_i + M\) for \(1<
      i+1<j-1\) such that \(a_{i+1}\equiv K_{i+1}\pmod\delta \).  Now pick
      \(a_{j-1,0}\ge a_{j-2}+M\) with \(a_{j-1,0}\equiv K_{j-1}\pmod\delta \),
      and inductively choose \(a_{j-1,k}\) and \(a_{j,k}\) for \(k<N\) to
      satisfy:
      \begin{itemize}
        \item
          \(a_{j,k}\ge a_{j-1,k}+M\)
        \item
          \(a_{j,k}\equiv K_j\pmod\delta \)
        \item
          \(a_{j-1,k+1}\ge a_{j,k}+M\)
        \item
          \(a_{j-1,k+1}\equiv K_{j-1}\pmod\delta \)
      \end{itemize}
      Now pick \(a_{j+1}\ge a_{j,N-1}+M\) with \(a_{j+1}\equiv K_{j+1}\pmod
      \delta \), and proceed inductively to pick \(a_{i+1}\ge a_i+M\) with
      \(a_{i+1}\equiv K_{i+1}\pmod\delta \) for \(j+1<i\le n\).
      Pictorially:
      \begin{tcd}
        a_1 \ar{r} & \cdots \ar{r} & a_{j-2} \ar{r}
        & a_{j-1,0} \ar{r} & a_{j,0}\ar[swap]{dl} &
        a_{j+1}\ar{r} & \cdots \ar{r} & a_n \\
        &&& a_{j-1,1} \ar{r} & a_{j,1}\ar[swap]{dl} \\
        &&& \vdots & \vdots \ar[swap]{dl} \\
        &&& a_{j-1,N-1} \ar{r} & a_{j,N-1} \ar{uuur}
      \end{tcd}
      where an arrow in the diagram indicates that the target is at least the
      source plus \(M\).

      For convenience we let \(b_k = (a_1,\ldots ,a_{j-2},a_{j-1,k})\) and
      \(c_\ell  = (a_{j,\ell },a_{j+1},\ldots ,a_n)\).
      Note now that for any \(k,\ell <N\) we have \(\models \chi (b_k,c_\ell
      )\). Furthermore if \(k\le \ell \) then \((b_k,c_\ell )\) satisfies
      \(x_1<\cdots < x_n\), so
      \(\models \psi_f (b_k,c_\ell )\)
      and thus \(\models \phi (b_k,c_\ell )\). Finally if \(k>\ell \) then
      \((b_k,c_\ell )\) satisfies \(x_{(j-1\; j)(1)}<\cdots < x_{(j-1\;
      j)(n)}\), so \(\models \psi _{f'}(b_k,c_\ell )\), and thus \(\models
      \neg\phi (b_k,c_\ell)\).

      Thus \(\models \phi (b_k,c_\ell )\) if and only if \(k\le \ell \), and we
      have constructed a ladder of size \(N\) for \(\phi \). So \(\phi \) has
      the order property and is thus unstable with respect to this partition of
      the variables, a contradiction. So no such \(\overline{a}\) exists, and
      we can compose \(\sigma \) with a transposition of adjacent elements and
      remain in \(\phi \). But such transpositions generate all of \(S_n\); so
      we may omit the ordering altogether and remain in \(\phi \), and thus
      \(\models \chi \rightarrow \phi \), as desired.
  \end{caselist}
  So we may replace \( \psi_f \) in the disjunction with a weaker
  \(L_\delta  \)-formula. Doing this to all disjuncts, we have written
  \(\phi \) as an \(L_\delta  \)-formula.
\end{rpf}
The remainder of our proof will make use of the \(F\)-sets of \cite{moosa04};
we briefly recall them here in the context of \(\mathbb{Z}^m\) where \(F\colon
\mathbb{Z}^m\to\mathbb{Z}^m\) is multiplication by \(d\).
\begin{definition}
  \label{def:fsets}
  A \emph{groupless \(F\)-set} in \(\mathbb{Z}^m\) is a finite union of
  translates of sums of sets of the form \(C(a;\delta )\) as defined above for
  \(a\in\mathbb{Z}^m\) and \(\delta >0\). An \emph{\(F\)-set} in
  \(\mathbb{Z}^m\) is a finite union of
  \[b + \sum_{i<n} C(a_i;\delta _i) + H
  \]
  for some \(b,a_i\in\mathbb{Z}^m\), \(\delta _i>0\), and \(H\le\mathbb{Z}^m\).
  The \emph{\(F\)-structure on \(\mathbb{Z}\)}, denoted
  \((\mathbb{Z},\mathcal{F})\), is the structure with domain \(\mathbb{Z}\) and
  a predicate for every \(F\)-set in every \(\mathbb{Z}^m\). (Note in
  particular that the graph of addition is a subgroup of \(\mathbb{Z}^3\), and
  hence an \(F\)-set; so \((\mathbb{Z},\mathcal{F})\) expands
  \((\mathbb{Z},+)\).)
\end{definition}
We now describe a simplification of \(d\)-sparsity that we will use to connect
stable sparse sets to \(F\)-sets.
\begin{lemma}
  \label{lemma:sparseincycles}
  Any \(d\)-sparse subset of \(\mathbb{Z}^m\) can be written as a finite union
  of translates of sets of the form
  \[\set{[\sigma _1^{e_1}] + \cdots + [\sigma _n^{e_n}] :
    e_1\le\cdots \le e_n} 
  \]
  where \(\sigma _i \in(\mathbb{Z}^m)^*\) all have the same length.
\end{lemma}
To see how this relates to \(F\)-sets, recall from \cref{rem:fsets} that
\([\sigma ^*] = C(\sigma ;\abs\sigma ) \cup\set0\); so a set of the above form
is (ignoring for the moment the case where some \(e_i = 0\)) a subset of the
groupless \(F\)-set \(C(\sigma _1;\abs{\sigma _1}) + \cdots +C (\sigma
_n;\abs{\sigma _n})\) that is cut out by some kind of order relation.
\begin{rpf}
  By \cref{fact:sparsechar} we can write \(A\) as a union of sets of the form
  \([u_0v_1^*\cdots v_n^*u_n]\) for \(u_i,v_i\in(\Sigma _\pm^m)^*\). Note
  first that we may assume all \(v_i\) across the union have the same length
  \(N\).  Indeed, let \(N\) be the least common multiple of the lengths of all
  the \(v_i\) across the union. We can then rewrite any \(v_i^*\) as
  \[\bigcup _{j< \ell } v_i^j(v_i^\ell )^*
  \]
  where \(\ell = \frac N{\abs{v_i}}\) (so \(\abs{v_i^\ell } = N\)). Using this
  to replace each \(v_i^*\) in the union and then distributing yields the
  desired expression for \(A\).

  It then suffices to show that given \(A = [u_0v_1^*\cdots v_n^*u_n]\) with
  each \(\abs{v_i} = N\) we can write \(A\) as a translate of a set of the form
  \(\set{[\sigma _1^{e_1}] + \cdots + [\sigma _n^{e_n}] : e_1\le\cdots \le e_n}
  \) with each \(\abs{\sigma _i} = N\).
  \begin{claim}
    We can write such \([u_0v_1^*\cdots v_n^*u_n]\) in the form \([a\tau
    _1^*\cdots \tau _n^*]\) where \(a\in(\mathbb{Z}^m)^*\) and each \(\tau
    _i\in (\mathbb{Z}^m)^*\) has length \(N\).
  \end{claim}
  \begin{rpf}
    We apply induction on \(n\); the base case is trivial.  For the induction
    step, use the induction hypothesis to write \([u_1v_2^*u_2\cdots
    u_{n-1}v_n^*u_n]\) as \([b\tau _2^*\cdots \tau _n^*]\).  Let \(x = [v_1] +
    d^{\abs{v_1}}[b] - [b]\); then the following is a telescoping sum:
    \[
      \underbrace{[u_0]+d^{\abs{u_0}}[b] +
      d^{\abs{u_0}}x}_y + \cdots + d^{\abs{u_0} + (k-1)\abs{v_1}}x
      =
      [u_0]+ d^{\abs{u_0}}[v_1] + \cdots  + d^{\abs{u_0} + (k-1)\abs{v_1}}[v_1]
      + d^{\abs{u_0} + k\abs{v_1}}[b]
    \]
    (The above equation is taken from a draft of \cite{bell19}; it was removed
    from the final paper.) Then if we let
    \begin{ea}
      a &=& y\cdot 0^{\abs{u_0}+\abs{v_1}-1} \\
      \tau _1 &=& x\cdot 0^{\abs{v_1}-1}
    \end{ea}
    (i.e.\ strings in \((\mathbb{Z}^m)^*\) whose first entries are \(y\) and
    \(x\) and whose later entries are the zero tuple) then
    \[[u_0v_1^kb] =
      \begin{cases}
        [a\tau _1^{k-1}]&\tif k>0 \\
        [u_0b]&\telse
      \end{cases}
    \]
    Hence if \(k>0\) then
    \[[u_0v_1^kbw]
      = [u_0v_1^kb] + d^{\abs{u_0v_1^kb}}[w]
      = [a\tau _1^{k-1}] + d^{\abs{a\tau _1^{k-1}}}[T_{\abs b}w]
      = [a\tau _1^{k-1}(T_{\abs b}w)]
    \]
    where \(T_i\sigma \) is the word obtained be replacing each letter
    \(\ell\in\mathbb{Z}^m \) appearing in \(\sigma \) with \(d^i\ell \). Hence
    \begin{ea}
      \relax[u_0v_1^*\cdots v_n^*u_n] &=& \set{[u_0bw]:w\in a_2^*\cdots
      a_n^*}\cup \set{[u_0v_1^kbw] : k\ge 1, w\in \tau _2^*\cdots \tau _n^*} \\
      &=& [u_0b\tau _2^*\cdots \tau _n^*] \cup \set{[a\tau _1^{k-1} (T_{\abs
      b}w)] : k\ge 1, w\in \tau _2^*\cdots \tau _n^*} \\
      &=& [u_0b\tau _2^*\cdots \tau _n^*] \cup [a\tau _1^*(T_{\abs b}\tau
      _2)^*\cdots (T_{\abs b}\tau _n)^*]
    \end{ea}
    And \(\abs{\tau _1} = \abs{v_1} = \abs{T_{\abs b}\tau _i} = N\) for all
    \(i\), as desired.
  \end{rpf}
  Note that given a set of the form \([a\tau _1^*\cdots \tau _n^*]\) with each
  \(\abs{\tau _i} = N\) we can rewrite it as \([a] + [(T_{\abs a}\tau
  _1)^*\cdots (T_{\abs a}\tau _n)^*] \). It then suffices to show that a set of
  the form \([\tau _1^*\cdots \tau _n^*]\) where each \(\tau
  _i\in(\mathbb{Z}^m)^*\) has length \(N\) can be written in the form
  \[\set{[\sigma _1^{e_1}] + \cdots + [\sigma _n^{e_n}] :
    e_1\le\cdots \le e_n} 
  \]
  with each \(\sigma _i\in(\mathbb{Z}^m)^*\) of length \(N\). For \(1\le i\le
  n\) let \(\sigma _i\in (\mathbb{Z}^m)^*\) be any string of length \(N\) such
  that \([\sigma _i] = [\tau _i] - \sum_{j=i+1}^n[\sigma _j]\). Then if
  \(e_1\le\cdots \le e_n\) then
  \begin{ea}
    \relax[\sigma _1^{e_1}] + \cdots + [\sigma _n^{e_n}]
    &=& [\sigma _1^{e_1}] + \\
    && [\sigma _2^{e_1}] + d^{Ne_1}[\sigma _2^{e_2-e_1}] + \\
    && \vdots \\
    && [\sigma _n^{e_1}] + d^{Ne_1}[\sigma _n^{e_2-e_1}] + \cdots +
    d^{Ne_{n-1}}[\sigma _n^{e_n-e_{n-1}}] \\
    &=& [\tau _1^{e_1}] + d^{Ne_1}[\tau _2^{e_2-e_1}] + \cdots +
    d^{Ne_{n-1}}[\tau _n^{e_n-e_{n-1}}] \\
    &=& [\tau _1^{e_1}\tau _2^{e_2-e_1}\cdots \tau _n^{e_n-e_{n-1}}]
  \end{ea}
  So \([\tau _1^*\cdots \tau _n^*] = \set{[\sigma _1^{e_1}] + \cdots + [\sigma
  _n^{e_n}] : e_1\le\cdots \le e_n}\), as desired.
\end{rpf}
The promised connection between stable sparse sets and \(F\)-sets:
\begin{lemma}
  Suppose \(A\subseteq \mathbb{Z}^m\) is \(d\)-sparse and stable in
  \((\mathbb{Z}^m,+)\).  Then \(A\) is definable in
  \((\mathbb{Z},\mathcal{F})\).
\end{lemma}
\begin{rpf}
  By \cref{lemma:sparseincycles} we can write \(A\) as a finite union of sets
  of the form
  \[\alpha +\set{[\sigma  _1^{e_1}] + \cdots +[\sigma _n^{e_n}] :
    e_1\le\cdots \le e_n}
  \]
  with \(\alpha \in\mathbb{Z}^m\) and each \(\sigma _i\in (\mathbb{Z}^m)^*\)
  has the same length \(N\). Take one such component of the union; we will show
  that it is contained in a set definable in \((\mathbb{Z},\mathcal{F})\) that
  is itself contained in \(A\), and hence can be replaced without changing the
  union.

  Since \(A\) is stable in \((\mathbb{Z}^m,+)\) and addition is associative we
  get that \(x_0 + x_1 + \cdots + x_n\in A\) is stable under any partition of
  the variables; thus so too is
  \[(x_1+\cdots +x_n\in A-\alpha )\wedge \bigwedge _{i=1}^n (x_i\in [\sigma
    _i^*])
  \]
  Thus \(X :=\set{(e_1,\ldots ,e_n) \in\mathbb{N}^n: [\sigma _1^{e_1}] + \cdots
  + [\sigma _n^{e_n}] \in A-\alpha }\) is a stable relation on \(\mathbb{N}^n\)
  under any partition of the variables. Furthermore if \(f\in S_n\) and we let
  \(\tau _i \in(\mathbb{Z}^m)^*\) be of length \(N\) such that \([\tau _i]=
  \sum_{j=i}^n [\sigma _{f(i)}]\), then as in the proof of
  \cref{lemma:sparseincycles} we get for \(e_{f (1)}\le\cdots \le e_{f (n)}\)
  we have
  \[(e_1,\ldots ,e_n)\in X\iff [\tau _1^{e_{f (1)}}\tau _2^{e_{f (2)}-e_{f
    (1)}}\cdots \tau _n^{e_{f (n)}-e_{f (n-1)}}] \in A-\alpha 
  \]
  Let \(\Lambda \supseteq \Sigma _\pm^m\) be any alphabet containing
  all the entries of the \(\tau _i\); so by \cref{rem:autorobust} we get that
  \(\set{\mu \in\Lambda ^* : [\mu ]\in A-\alpha }\) is regular.  So, by the
  proof of \cref{prop:sparseintersregpresb} we get that \([\tau _1^{t_1}\tau
  _2^{t_2}\cdots \tau _n^{t_n}]\in A-\alpha \) can be expressed by a Boolean
  combination of congruences and equalities between a \(t_i\) and a constant.
  But a congruence or equality between \(e_{f (i+1)}-e_{f (i)}\) and a constant
  \(k\) can be expressed as a congruence or equality between \(e_{f (i+1)}\)
  and \(S^k(e_{f (i)})\), and is thus expressible by an \(L_\delta \)-formula
  for some \(\delta \); furthermore by taking disjunctions and LCMs we may
  assume all congruences that occur have the same modulus \(\delta \).  So
  \[[\tau _1^{e_{f (1)}}\tau _2^{e_{f (2)}-e_{f
    (1)}}\cdots \tau _n^{e_{f (n)}-e_{f (n-1)}}] \in A-\alpha 
  \]
  can be expressed as an \(L_\delta \)-formula for some \(\delta \), and hence
  so too can \((e_1,\ldots ,e_n)\in X\) as long as \(e_{f (1)}\le\cdots \le
  e_{f (n)}\). So taking disjunctions over possible orderings of the \(e_i\),
  (and LCMs of the resulting \(\delta \)) we see that \(X\) can be defined by a
  quantifier-free \(L_{\delta ,<}\)-formula that is stable under any partition
  of the variables. So by previous proposition we get that \(X\) can be defined
  by a quantifier-free \(L_\delta \)-formula.

  Let \(\mathbbm1 \in \mathbb{Z}^m\) be the tuple all of whose entries are
  \(1\).  I claim that the map \(\mathbb{N}\to \mathbb{Z}^m \) given by
  \(e\mapsto [\mathbbm1^{Ne}]\) defines an interpretation of \((\mathbb{N},0,S,
  P_\delta )\) in \((\mathbb{Z} ,\mathcal{F})\).  Indeed, the image is
  definable: it's simply \(C([\mathbbm1^N];N)\cup\set0\). Furthermore the
  unnested atomic \(L_\delta \)-formulas all map to definable sets in
  \((\mathbb{Z} ,\mathcal{F})\):
  \begin{itemize}
    \item
      \(P_\delta (x)\) can be expressed as a condition of \([\mathbbm1^{Nx}]\)
      by demanding that it lie in \(C([\mathbbm1^{N\delta }]; N\delta
      )\cup\set0 \).
    \item
      \(y = Sx\) can be expressed by requiring that
      \(\col{\relax[\mathbbm1^{Nx}]}{\relax[\mathbbm1^{Ny}]}\in
      \pars*{C\pars*{\col{\relax[\mathbbm1^N]}{d^N[\mathbbm1^N]}; N}
      \cup\set0}+ \col 0{\relax[\mathbbm1^N]}\).
  \end{itemize}
  Furthermore the map \([\mathbbm1^{Ne}]\mapsto [\sigma _i^e]\) is definable in
  \((\mathbb{Z} ,\mathcal{F})\) for each \(i\): its graph is simply
  \(C\pars*{\col {\relax[\mathbbm1^N]}{\relax[\sigma _i]}; N}\cup\set0\).
  (Recall that \(\abs{\sigma _i} = N\).) Then since \(X\) is definable in
  \((\mathbb{N},0,S,P_\delta )\) we get that
  \[Y := \set*{\sum_{i=1}^n[\sigma _i^{e_i}] : (e_1,\ldots ,e_n)\in X}
    = (A-\alpha )\cap ([\sigma _1^*]+\cdots +[\sigma _n^*])
  \]
  is definable in \((\mathbb{Z} ,\mathcal{F})\). (Recall that addition is
  definable in \((\mathbb{Z},\mathcal{F})\).) But \([\sigma _1^{e_1}]+\cdots
  +[\sigma _n^{e_n}]\in A-\alpha \) if \(e_1\le\cdots \le e_n\); so
  \[\alpha +\set{[\sigma _1^{e_1}] + \cdots + [\sigma _n^{e_n}] : e_1\le\cdots
    \le e_n}\subseteq \alpha+Y \subseteq A
  \]
  So we can replace \(\alpha+\set{[\sigma _1^{e_1}] + \cdots + [\sigma
  _n^{e_n}] : e_1\le\cdots \le e_n }\) in the union defining \(A\) with
  \(\alpha+Y \).  Doing this for all such terms in the union, we can write
  \(A\) as a union of sets definable in \((\mathbb{Z} ,\mathcal{F})\); so \(A\)
  is definable in \((\mathbb{Z} ,\mathcal{F})\).
\end{rpf}
Our theorem now follows:
\begin{lrpf}{thm:charstablesparse}
  By previous lemma \(A\) is definable in \((\mathbb{Z},\mathcal{F})\).  But by
  \cite[Theorem A]{moosa04} \((\mathbb{Z},\mathcal{F})\) admits quantifier
  elimination. So \(A\) is definable by a Boolean combination of \(F\)-sets,
  say in disjunctive normal form; we must show the \(F\)-sets can be taken to
  be groupless. Take one disjunct
  \[\bigcap _{i<k} B_i \setminus \bigcup _{j<\ell } C_j
  \]
  where the \(B_i,C_j\) are \(F\)-sets. By \cref{lemma:sparseincycles} \(A\) is
  contained in a finite union of translates of finite sums of \(C(\sigma
  _i;\delta _i)\), and hence in a groupless \(F\)-set \(\widehat{A}\). So if
  \(k>0\) we may replace every \(B_i\) and \(C_j\) in our disjunct with
  \(B_i\cap\widehat{A}\) and \(C_j\cap\widehat{A}\), respectively, and the
  result of the disjunction will still be \(A\). If \(k = 0\) we instead
  replace our disjunct with \(\widehat{A}\setminus \bigcup _{j<\ell }(C_j\cap
  \widehat{A})\), and again the result of the disjunction is still \(A\). But
  \(B_i\cap\widehat{A},C_j\cap\widehat{A}\) are intersections of \(F\)-sets,
  and hence themselves \(F\)-sets by \cite[Proposition 3.9]{moosa04}.
  Furthermore \(\widehat{A}\) is \(d\)-sparse by \cref{rem:fsets}; so
  \(B_i\cap\widehat{A},C_j\cap\widehat{A}\) cannot contain a translate of a
  subgroup, and hence are groupless \(F\)-sets.  Applying the above replacement
  to every disjunct, we get that \(A\) is a Boolean combination of groupless
  \(F\)-sets, i.e.\ translates of sums of \(C(a;\delta)\), as desired.
\end{lrpf}
We conclude by pointing out that combined with \cite{moosa04} we obtain the
following characterization of the stable \(d\)-sparse sets:
\begin{corollary}
  \label{cor:sparsestablechar}
  Suppose \(A\subseteq \mathbb{Z}^m\) is \(d\)-sparse. \Tfae{}:
  \begin{enumerate}
    \item
      \(\thy(\mathbb{Z},+,A)\) is stable.
    \item
      \(A\) is stable in \((\mathbb{Z}^m,+)\).
    \item
      \(A\) is a finite Boolean combination of translates of sums of sets of
      the form \(C(a;\delta )\).
    \item
      \(A\) is definable in \((\mathbb{Z},\mathcal{F})\).
  \end{enumerate}
\end{corollary}
\begin{rpf}
  (1)\(\implies \)(2) is clear, (2)\(\implies \)(3) is
  \cref{thm:charstablesparse}, (3)\(\implies \)(4) is clear, and (4)\(\implies
  \)(1) is by the fact (Theorem A of \cite{moosa04}) that
  \(\thy(\mathbb{Z},\mathcal{F})\) is stable.
\end{rpf}

\section{Beyond sparsity: the non-generic case}
\label{sec:nonsparse}
In the previous section we characterized the \(d\)-sparse sets that are stable
in \((\mathbb{Z}^m,+)\). So the question of which automatic sets are stable in
\((\mathbb{Z}^m,+)\) reduces to the non-sparse case. We begin to study this
problem in this section, restricting our attention to subsets of
\(\mathbb{Z}\).

As an example of a non-sparse automatic set that is stable in
\((\mathbb{Z},+)\), consider a coset of a subgroup, say \(A = r+s\mathbb{Z}\)
where \(s>0\). Then \(A\) is stable in \((\mathbb{Z},+)\) since it's definable
in \((\mathbb{Z},+)\).  It isn't \(d\)-sparse: the number of \(a\in A\) with
\(d^{-k}< a < d^k\) grows exponentially with \(k\), so the set of canonical
representations of \(A\) isn't sparse. It is \(d\)-automatic: see \cite[Theorem
5.4.2]{allouche03} (though recall as mentioned in \cref{rem:autorobust} that
they use a different convention for representing integers, so the automaton
will be slightly different).

One can also take Boolean combinations of cosets and the stable sparse sets of
the previous section to get further examples, as long as the result isn't
\(d\)-sparse. But all examples produced in this way will be ``generic'':
\begin{definition}
  We say \(A\subseteq \mathbb{Z}\) is \emph{generic} if some
  finite union of additive translates of \(A\) covers \(\mathbb{Z}\).
\end{definition}
We show that in the non-sparse setting all stable automatic sets are generic.
\begin{theorem}
  \label{thm:stablenongen}
  Suppose \(A\subseteq \mathbb{Z}\) is \(d\)-automatic and not \(d\)-sparse.
  If \(A\) is stable in \((\mathbb{Z},+)\) then \(A\) is generic.
\end{theorem}
It will be easier to work first in \(\mathbb{N}\), and in particular to use
\(\Sigma =\set{0,\ldots ,d-1}\) for our representations rather than \(\Sigma
_\pm =\set{-d+1,\ldots ,d-1}\); the main advantage to doing so is that whenever
\(\sigma ,\tau \in\Sigma ^*\) have the same length we have \(\sigma = \tau \iff
[\sigma ] = [\tau ]\). (Note that the same does not hold in \(\Sigma _\pm\):
for example \([(d-1)0] = [(-1)d]\).) Recall from \cref{rem:autorobust} that
\(A\subseteq \mathbb{N}\) is a \(d\)-automatic subset of \(\mathbb{Z}\) if and
only if it is a \(d\)-automatic subset of \(\mathbb{N}\) in the classical
sense; i.e.\ \(\set{\sigma \in\Sigma ^*: [\sigma ]\in A}\) is regular.  Note
also that if \(A\subseteq \mathbb{N}\) then the canonical representations of
the elements of \(A\) all lie in \(\Sigma ^*\), and up to trailing zeroes these
are the only representations over \(\Sigma \) of elements of \(A\). So
\(A\subseteq \mathbb{N}\) is \(d\)-sparse as a subset of \(\mathbb{Z}\) if and
only if \(\set{\sigma \in \Sigma ^*: [\sigma ]\in A, \sigma \text{ has no
trailing zeroes}}\) is sparse.  On the other hand stability and genericity when
relativized to \(\mathbb{N}\) give something new:
\begin{definition}
  We say \(A\subseteq \mathbb{N}\) is \emph{stable in \(\mathbb{N}\)} if
  \(x+y\in A\) is a stable relation on \(\mathbb{N}\). We say \(A\) is
  \emph{generic in \(\mathbb{N}\)} if some finite union of (possibly negative)
  translates of \(A\) covers \(\mathbb{N}\).
\end{definition}
We will first focus on proving:
\begin{proposition}
  \label{prop:stablenatnongen}
  Suppose \(A\subseteq \mathbb{N}\) is \(d\)-automatic and not \(d\)-sparse.
  If \(A\) is stable in \(\mathbb{N}\) then \(A\) is generic~in~\(\mathbb{N}\).
\end{proposition}
We begin with a characterization of the generic \(d\)-automatic sets.
\begin{lemma}
  \label{lemma:nongenforbsuf}
  Suppose \(A\subseteq \mathbb{N}\) is \(d\)-automatic; let \(L\subseteq \Sigma
  ^*\) be the set of representations of elements of \(A\). Then \tfae{}:
  \begin{enumerate}
    \item
      \(A\) is generic in \(\mathbb{N}\).
    \item
      For any \(r,s\in\mathbb{N}\), every \(\tau \in\Sigma ^*\) occurs as a
      suffix of a word in \(L\) of length \(r+sk\) for some \(k\ge0\).
  \end{enumerate}
\end{lemma}
In other words, \(A\) is not generic in \(\mathbb{N}\) if and only if there are
\(r,s\in\mathbb{N}\) such that \(L\cap\Sigma ^{(r+s\mathbb{N})}\) has a
forbidden suffix.
\begin{rpf}
  Note first that \(A\) is not generic in \(\mathbb{N}\) if and only if there
  are arbitrarily large gaps in \(A\) (i.e.\ runs of naturals not in \(A\)).
  \begin{description}
    \pip12
      Suppose we are given \(\tau ,r,s\) such that \(\tau \) is a forbidden
      suffix for \(L\cap\Sigma ^{(r+s\mathbb{N})}\). Then if \(r+sk>\abs\tau \)
      then \(A\) is disjoint from
      \[[\Sigma ^{(r+sk-\abs\tau) }\tau ] = \set{b\in\mathbb{N}:
      d^{r+sk-\abs\tau }[\tau] \le b< d^{r+sk-\abs\tau }([\tau ]+1)]}
      \]
      So \(A\) has a gap of size \(d^{r+sk-\abs\tau }\). So as \(k\to\infty \)
      we get arbitrarily large gaps in \(A\); so \(A\) isn't generic in
      \(\mathbb{N}\).
    \pip21
      Suppose \(A\) isn't generic in \(\mathbb{N}\).  Let \(\$\) be a letter
      not in \(\Sigma \); we will use \(\$\) as a separator. Consider the set
      \(S\subseteq (\Sigma \cup\set\$)^*\) of \(0^m\$\tau \) for \(m <\omega \)
      and \(\tau \in\Sigma ^*\) with the property that \(\Sigma ^{(m)}\tau \cap
      L=\emptyset \); in other words, if we replace each zero with any letter
      and delete the separator, the result is never in \(L\). So \(0^m\$\tau\in
      S\) if and only if \(\tau \) is a forbidden suffix for \(L\cap \Sigma
      ^{(m+\abs\tau) }\). Then \(S\) is regular: it's not too hard to construct
      a \emph{non-deterministic finite automaton} (NFA) for the complement,
      which suffices (see e.g.\ \cite[Section 2.2]{yu97}). Since there are
      arbitrarily large gaps in \(A\) we get that there are elements
      \(0^m\$\tau \in S\) with \(m\) arbitrarily large.  Indeed, suppose we are
      given \(m\). Find a gap of size \(2d^m\); then this gap will contain two
      multiples of \(d^m\), say \(a,a+d^m\). Then if \(\tau \in\Sigma ^*\) is
      such that \([\tau ] = \frac a{d^m}\) then \(\tau \) is a forbidden suffix
      for \(L\cap\Sigma ^{(m+\abs\tau) }\); so \(0^m\$\tau \in S\).

      Recall the pumping lemma for regular languages (see
      \cite[Lemma 4.1]{yu97}): if \(R\) is regular then there is a
      \emph{pumping length} \(p>0\) such that if \(\mu \in R\) has length \(\ge
      p\) then we can write \(\mu =uvw\) such that
      \begin{itemize}
        \item
          \(v\ne\epsilon \)
        \item
          \(\abs{uv}\le p\)
        \item
          \(uv^*w\subseteq R\).
      \end{itemize}

      Pick \(0^m\$\tau \in S\) with \(m\) bigger than the pumping length of
      \(S\).  Then by the pumping lemma we can write \(m = r+s\) so that
      \(0^r(0^s)^*\$\tau \subseteq S\); so \(\tau \) is a forbidden suffix
      for \(L\cap\Sigma ^{(r+\abs\tau +s\mathbb{N})}\).
      \qedhere
  \end{description}
\end{rpf}
The following technical lemma is the source of instability in
\cref{prop:stablenatnongen}. For \(K<\omega \) we define a partial binary
operation \(+_K\) on \(\Sigma ^*\) by setting \(\sigma +_K\tau \) to be the
unique representation of \([\sigma ]+[\tau ]\) of length \(K\), if one exists.
\begin{lemma}
  \label{lemma:nongenspecialcase}
  Suppose \(L\subseteq \Sigma ^*\) is regular but not sparse, and satisfies
  \(L=L^*\) and
  \[\label{cond:forbsuff}
    \text{there are \(r,s\in\mathbb{N}\) such that \(L\cap\Sigma
    ^{(r+s\mathbb{N})}\) is infinite and has a forbidden suffix \(\sigma
    \).} \tag{\dag}
  \]
  Then for all \(N<\omega \) there is \(K<\omega \) such that the binary
  relation \(x+_Ky\in L\) on \(\Sigma ^*\) has an \(N\)-ladder.
\end{lemma}
\begin{rpf}
  Pick \(\sigma ,r,s\) as in (\ref{cond:forbsuff}). Since \(L\cap\Sigma
  ^{(r+s\mathbb{N})}\) is infinite, there is \(a\in\Sigma ^{(\abs\sigma) }\)
  that occurs as a suffix of some element of \(L\cap\Sigma
  ^{(r+s\mathbb{N})}\).  Suppose \([a]\le[\sigma ]\); we will see at the end
  how to modify the argument in the case \([a]>[\sigma ]\).

  Pick such \(a\) maximal under \(\le_\mathbb{N}\) (the preorder induced from
  the ordering on \(\mathbb{N}\)); so if \(a'\in\Sigma ^{(\abs\sigma) }\) has
  \([a'] = [a]+1\) then \(a'\) does not occur as a suffix of some element of
  \(L\cap\Sigma ^{(r+s\mathbb{N})}\).  (Note such \(a'\) exists since \([a]<
  [\sigma ]<d^{\abs\sigma }\), and hence \([a]+1<d^{\abs\sigma }\) can be
  represented by a string of length \(\abs\sigma \).) Consider the set \(S\) of
  \(\tau \in L\cap\Sigma ^{(r+s\mathbb{N})}\) with \(a\) as a suffix such that
  \(\tau \) is \(\le_\mathbb{N}\)-maximal among the elements of \(L\) ending in
  \(a\) that are of the same length as \(\tau \).  Then \(S\) is infinite:
  since \(L = L^*\) and \(a\) occurs as a suffix of some \(\mu \in L\cap\Sigma
  ^{(r+s\mathbb{N})}\), we get that \(\mu ^{1+s\mathbb{N}}\subseteq L\cap\Sigma
  ^{(r+s\mathbb{N})}\) also has \(a\) as a suffix, and hence that \(S\)
  contains a word of length \((1+sk)\abs\mu \) for \(k<\omega \).  Furthermore
  \(S\) is regular: using the fact that \(\set{(\mu ,\nu)\in (\Sigma ^2)^* :
  [\mu ]\le[\nu ]}\) and \(\Sigma ^*a\) are regular, one can construct an NFA
  for the complement of \(S\).  So by the pumping lemma \(S\) contains a set of
  the form \(uv^*w\) with \(v\ne\epsilon \). By prepending a power of \(v\) to
  \(w\) we may assume \(\abs w\ge \abs a\), and in particular that \(w\) has
  \(a\) as a suffix (and is non-empty).

  Since \(L = L^*\) and \(uv^*w\subseteq S\subseteq L \) we get that
  \(L\supseteq (uv^*w)^*\supseteq u\set{wu, v}^*w\). This, together with the
  maximality of elements of \(S\), the fact that \(a'\) is a forbidden suffix
  for \(L\cap\Sigma ^{(r+s\mathbb{N})}\), and the fact that \(\abs{uv^*w}\in
  r+s\mathbb{N}\), will be enough to construct our ladder.

  Pick \(n,m\) such that \(n\abs{wu} = m\abs v\); then \(u(wu)^nw\in L\) and
  ends in \(a\), so since \(uv^mw\in S\) and \(\abs{uv^mw} = \abs{u(wu)^nw}\)
  we get that \([u(wu)^nw]\le[uv^mw]\), and hence that \([(wu)^n]\le [v^m]\).
  \begin{caselist}
    \case
      Suppose \([(wu)^n]<[v^m]\); then since \([v^m]-[(wu)^n]\le [v^m]\) there
      is \(\alpha \in \Sigma ^{(m\abs v)}\) such that \([\alpha ] =
      [v^m]-[(wu)^n]>0\). We let
      \begin{ea}
        K &=& \abs{u} + Nm\abs{v} + \abs w \in r+s\mathbb{N}\\
        d_i &=& u(wu)^{n(N-i)}v^{mi}w \\
        e_i &=& 0^{\abs u}\alpha ^{N-i}
      \end{ea}
      for \(i \le N\). Then \(d_i+_Ke_j\) is defined for all \(i,j\); i.e.\
      \([d_i]+ [e_j]\) has a representation of length \(K\). Indeed,
      \(\abs{e_j}\le K-\abs w\le K-\abs a\); so \([e_j]< d^{K-\abs a}\). So if
      we write \(d_i = \tau a\) for some \(\tau \) (possible since \(d_i\) has
      \(w\), and hence \(a\), as a suffix) then
      \[[d_i] + [e_j] < [\tau a]+d^{K-\abs a} = [\tau a] + d^{\abs\tau }
        = [\tau a'] < d^K
      \]
      since \(\abs{\tau a'} = \abs{d_i} = K\). So \([d_i] + [e_j]\) has a
      representation of length \(K\), and \(d_i+_Ke_j\) is defined. In fact the
      above proof shows that \(d_i+_Ke_j\) has either \(a\) or \(a'\) as a
      suffix.

      Since \([\alpha ]>0\) it is clear that the \(e_i\) are strictly
      decreasing. Suppose \(i>j\); then \([d_i+_Ke_j] = [d_i] + [e_j] > [d_i] +
      [e_i] = [uv^{mN}w]\). So if \(d_i+_Ke_j\) has \(a\) as a suffix then
      since \(uv^{mN}w\in S\) and \(d_i+_Ke_j\) has the same length, has \(a\)
      as a suffix, and represents a strictly larger number, we get that
      \(d_i+_Ke_j\notin L\). Otherwise as noted above we get that \(d_i+_K
      e_j\) has \(a'\) as a suffix, in which case \(d_i+_Ke_j\notin L\) since
      \(a'\) is a forbidden suffix for \(L\cap\Sigma ^{(r+s\mathbb{N})}\) and
      \(\abs{d_i+_Ke_j} = K\in r+s\mathbb{N}\). Conversely suppose \(i\le j\);
      then \(d_i +_K e_j = uv^{m(N-j)}(wu)^{n(j-i)}v^{mi}w\in
      u\set{wu,v}^*w\subseteq L\). So the \(d_i,e_i\) form an \(N\)-ladder for
      \(x+_Ky\in L\).
    \case
      Suppose \([(wu)^n] = [v^m]\); so \((wu)^n = v^m\). Then \( uv^*w
      \supseteq u((wu)^n)^*w = ((uw)^n)^*uw\); so if we let \(u'= \epsilon \),
      \(v' = (uw)^n\), and \(w' = uw\) then \( u'(v')^*w'\subseteq
      uv^*w\subseteq S\). Furthermore \(v'\ne\epsilon \) since \(\abs{v'} =
      n\abs{uw} = m\abs{v} > 0\), and \(v'\in L\) since \(L=L^*\) and \(uw\in
      uv^*w\subseteq L\). So we may replace \(u,v,w\) with \(u',v',w'\)
      respectively, and we may thus assume that \(u=\epsilon \) and \(v,w\in
      L\).  (Recall that the only requirement we had of \(u,v,w\) was that
      \(uv^*w\subseteq S\) and \(v\ne\epsilon \).)

      By \cite[Proposition 7.1]{bell19} since \(L\) isn't sparse there are
      \(x,y_1,y_2,z\in\Sigma ^*\) with \(y_1,y_2\) distinct, non-trivial, and
      of the same length such that \(x\set{y_1,y_2}^*z\subseteq L\). Let
      \(b=xy_1z\) and \(c=xy_2z\); so \(\abs b= \abs c\) with \(b,c\in L\) and
      \(b\ne c\).  By replacing \(b,c,v\) with powers thereof we may assume
      \(\abs b = \abs c = \abs v\). Then since \(b\ne c\) we get that one of
      \(b,c\), without loss of generality say \(b\), has \(b\ne v\), and thus
      \([b]\ne[v]\). Note since \(L=L^*\) that \(L\supseteq \set{b,v}^*w\).

      Since \(vw\in S\) and since \(bw\) has the same length as \(vw\), has
      \(a\) as a suffix, and lies in \(L\), we get that \([bw]\le [vw]\). So
      \([b]\le [v]\), and since \(b\ne v\) we get \([b]<[v]\). Then since
      \([v]-[b]<[v]\) there is \(\alpha \in\Sigma^{(\abs{v})}\) such that
      \([\alpha ] = [v]-[b]\). We then let
      \begin{ea}
        K &=& N\abs v+\abs w \in r+s\mathbb{N}\\
        d_i &=& b^{N-i}v^iw  \\
        e_i &=& \alpha ^{N-i}
      \end{ea}
      for \(i\le N\).  Then by an argument identical to the previous case the
      \(d_i,e_i\) form an \(N\)-ladder for \(x+_Ky\in L\).
  \end{caselist}
  The case \([a]>[\sigma ]\) is similar; we outline it here. We take minimal
  such \(a\) under \(\le_\mathbb{N}\), and define \(S\) to be the set of
  \(\tau\in L\cap\Sigma ^{(r+s\mathbb{N})} \) ending in \(a\) that are
  \(\le_\mathbb{N}\)-minimal among the elements of \(L\) ending in \(a\) that
  are of the same length as \(\tau \).  Then \(S\) is again infinite and
  regular, and thus contains a set of the form \(uv^*w\); we again assume \(w\)
  has \(a\) as a suffix. If \(n\abs{wu}=m\abs v\) then dually to before we get
  \([(wu)^m]\ge[v^n]\). If \([(wu)^n]>[v^m]\), say with \(\alpha \in \Sigma
  ^{(m\abs v)}\) with \([\alpha ] = [(wu)^n] - [v^m]>0\), then we'd like to let
  \begin{ea}
    K &=& \abs{u} + Nm\abs{v}+\abs w \\
    d_i &=& u(wu)^{n(N-i)}v^{mi}w \\
    e_i &=& 0^{\abs u}(-\alpha) ^{N-i}
  \end{ea}
  and claim this as our ladder.  Unfortunately we're working over \(\Sigma \),
  not \(\Sigma _\pm\), so we can't allow the \(e_i\) to use negative digits.
  This is easily fixed, however: note for all \(i,j\) that \([d_i] \ge d^{\abs
  u+Nm\abs v}\ge -[e_j]\) (since \([a]\ne 0\) and \(w\), and hence \(d_i\), has
  \(a\) as a suffix). So we can take \(d_i',e_i'\in\Sigma ^*\) such that
  \([d_i'] = [d_i] - d^{\abs u+Nm\abs v}\) and \([e_i'] = e_i + d^{\abs
  u+Nm\abs v}\).  Then \([d_i'] + [e_j'] = [d_i] + [e_j]\), and now as before
  one can show that \(d_i'+_Ke_j'\) is always defined and is in \(L\) if and
  only if \(i\le j\).

  If \([(wu)^n] = [v^m]\) we do a similar trick. As before we may assume
  \(u=\epsilon \) and \(v,w\in L\), and we get some \(b\in L\) with \(\abs b =
  \abs v\) and \([b]\ne [v]\); dually to before we get \([b]>[v]\), say with
  \(\alpha \in\Sigma ^{(\abs v)}\) such that
  \([\alpha ] = [b]-[v]\). Our initial attempt at a ladder will now be:
  \begin{ea}
    K &=& N\abs{v} + \abs w \\
    d_i &=& b^{N-i}v^iw \\
    e_i &=& (-\alpha )^{N-i}
  \end{ea}
  Now we have \([d_i] \ge d^{N\abs v} \ge -[e_j]\); so we can pull the same
  trick to turn the \(d_i,e_i\) into a ladder.
\end{rpf}
% We still need to show that we can reduce the case of a general non-generic,
% \(d\)-automatic, and non-\(d\)-sparse \(A\subseteq \mathbb{N}\) to the above
% lemma. We do so by taking the automaton for \(A\) and examining a related
% automaton that will recognize a language satisfying the hypotheses of
% \Cref{lemma:nongenspecialcase}.
Suppose \(M= (Q,q_0,F,\delta )\) is a DFA over \(\Sigma \). For \(q\in Q\) we
let \(L_q = \set{\sigma \in\Sigma ^* : \delta (q,\sigma ) = q}\); that is,
\(L_q\) is the set of words which take state \(q\) back to state \(q\) in
\(M\). Note that \(L_q\) is regular: it is recognized by the automaton
\((Q,q,\set{q},\delta )\).
\begin{lemma}
  \label{lemma:stateconditions}
  Suppose \(A\subseteq \mathbb{N}\) is \(d\)-automatic but not \(d\)-sparse;
  suppose \(A\) is not generic in \(\mathbb{N}\). Fix an automaton \(M=
  (Q,q_0,F,\delta )\) that recognizes the set of representations over \(\Sigma
  \) of elements of \(A\).  Then there is a non-dead \(q\in Q\) such that
  \(L_q\) satisfies the hypotheses of \cref{lemma:nongenspecialcase}: namely
  \(L_q\) is regular but not sparse, \(L_q = L_q^*\), and \(L_q\) satisfies
  (\ref{cond:forbsuff}).
\end{lemma}
(We say \(q\in Q\) is a \emph{dead state} if there is no \(\sigma \) such that
\(\delta (q,\sigma )\in F\).)
\begin{rpf}
  Note we always have that \(L_q\) is regular and \(L_q = L_q^*\); so we only
  need non-sparsity and (\ref{cond:forbsuff}).  We first note some facts about
  how non-sparsity and (\ref{cond:forbsuff}) interact with the \(L_q\).
  \begin{claim}
    \label{claim:lqprops}
    \begin{menumerate}
      \item
        If \(q\) is a finish state of \(M\) and \(L_q\) is infinite then
        \(L_q\) satisfies (\ref{cond:forbsuff}).
      \item
        There is a non-dead \(q\) such that \(L_q\) isn't sparse.
      \item
        If \(q,q'\) are states in \(M\) with a path from \(q\) to \(q'\) and
        vice-versa then \(L_q\) is sparse if and only if \(L_{q'}\) is.
    \end{menumerate}
  \end{claim}
  \begin{rpf}
    \begin{menumerate}
      \item
        Let \(L\) be the set of representations of elements of \(A\), and fix
        \(\mu \in\Sigma ^*\) such that \(\delta (q_0,\mu ) = q\). (We may
        assume such \(\mu \) exists: otherwise we can remove \(q\) from \(M\)
        without changing the set recognized by \(M\).) By non-genericity of
        \(A\) in \(\mathbb{N}\) and \cref{lemma:nongenforbsuf} there is some
        forbidden suffix for \(L\cap\Sigma ^{(r+s\mathbb{N})}\).  Note that if
        \(\tau \) is a forbidden suffix for \(L\cap \Sigma ^{(r+s\mathbb{N})}\)
        then \(\tau 0^t\) is a forbidden suffix for \(L\cap\Sigma
        ^{(r+t+s\mathbb{N})}\) (since \(L\) is closed under removing trailing
        zeroes). So there is a forbidden suffix for \(L\cap \Sigma
        ^{(r'+s\mathbb{N})}\) for any \(r'\ge r\); pick \(r'\) such that
        \(L_q\cap \Sigma ^{(r'+s\mathbb{N}-\abs\mu) }\) is infinite.  Then
        since \(q\) is a finish state the forbidden suffix for \(L\cap\Sigma
        ^{(r'+s\mathbb{N})}\) is also a forbidden suffix for \(L_q\cap \Sigma
        ^{(r'+s\mathbb{N}-\abs\mu) }\). So \(L_q\) satisfies
        (\ref{cond:forbsuff}).
      \item
        By \cite[Proposition 7.1]{bell19} there is a non-\(M\)-dead state \(q\)
        and distinct non-empty \(u,v\in\Sigma ^*\) such that \(\delta (q,u) =
        \delta (q,v) = q\) and \(\delta (q,x)\ne q\) for \(x\) any proper
        non-empty prefix of \(u\) or \(v\). Then taking \(b,c\) to be powers of
        \(u,v\) respectively such that \(\abs b=\abs c\), we get that \(b\ne
        c\) (otherwise \(u\) or \(v\) would be a prefix of the other); also
        \(\delta (q,a) = \delta (q,b) = q\), so \(b,c\in L_q\). So
        \(L_q\supseteq \set{b,c}^*\), and a quick computation shows that
        \(L_q\) isn't sparse.
      \item
        Let \(\delta (q,\mu )= q'\) and \(\delta (q',\nu ) = q\). Suppose
        \(L_q\) isn't sparse. Then \(L_{q'}\supseteq \nu L_q\mu \) also isn't
        sparse.
        \qedhere
    \end{menumerate}
  \end{rpf}
  By \cref{claim:lqprops} (2) there is \(q\) such that \(L_q\) isn't sparse and
  \(a\) has a path to a finish state \(q'\). If \(L_{q'}\) isn't sparse then by
  \cref{claim:lqprops} (1) we're done; suppose then that it is sparse.  We show
  in this case that there is a forbidden infix for \(L_q\) (i.e.\ some \(\sigma
  \) that does not appear as a substring of any element of \(L_q\)), and hence
  in particular that \(L_q\) satisfies (\ref{cond:forbsuff}) with \(r=0\) and
  \(s=1\).

  Note that there is no path from \(q'\) to \(q\), else by \cref{claim:lqprops}
  (3) \(L_{q'}\) wouldn't be sparse.  Enumerate the states of \(M\) with a path
  to \(q\) (and hence to \(q'\)) as \((q_i : i <n)\). Inductively pick \(\sigma
  _i\in\Sigma ^*\) as follows: if \(\delta (q_i,\sigma _0\cdots \sigma
  _{i-1})\) has no path to \(q\) we let \(\sigma _i =\epsilon \), and
  otherwise we pick \(\sigma _i\) such that \(\delta (q_i,\sigma _0\cdots
  \sigma _i) = q'\). Note then that \(\delta (q_i,\sigma _0\cdots \sigma _i)\)
  has no path to \(q\); hence neither does \(\delta (q_i,\sigma _0\cdots \sigma
  _{n-1})\). Let \(\tau = \sigma _0\cdots \sigma _{n-1}\). We have shown that
  if \(r\) is a state with a path to \(q\) (so one of the \(q_i\)) then
  \(\delta (r,\tau ) \) has no path to \(q\). Clearly if \(r\) has no path to
  \(q\) then neither does \(\delta (r,\tau )\). Hence for all \(r\in Q\) we get
  that \(\delta (r,\tau)\) has no path to \(q\); that is, \(\tau \) is a
  forbidden infix.
\end{rpf}
\begin{lrpf}{prop:stablenatnongen}
  Suppose \(A\subseteq \mathbb{N}\) is \(d\)-automatic and neither \(d\)-sparse
  nor generic in \(\mathbb{N}\).  Fix a minimal automaton \(M = (Q,
  q_0,F,\delta )\) for the set of representations over \(\Sigma \) of elements
  of \(A\). (The \emph{minimal automaton} of a regular language \(L\) is an
  automaton recognizing \(L\) where all states are reachable from the start
  state and such that given distinct \(q,q'\in Q\) there is \(\nu \) such that
  \(\delta (q,\nu )\in F\) if and only if \(\delta (q',\nu )\notin F\). Such
  automata exist and are unique: see the proof of the right-to-left direction
  of \cite[Theorem 4.7]{yu97}.)

  By \cref{lemma:stateconditions} there is a non-dead \(q\) such that
  \(L_q\) satisfies the hypotheses of \cref{lemma:nongenspecialcase}.  Using
  minimality,
  for each \(q'\ne q\) let \(\sigma _{q'}\in\Sigma
  ^*\) and \(\epsilon _{q'}\in \set{0,1}\) be such that \((\delta (q,\sigma
  _{q'})\in F)^{\epsilon _{q'}}\wedge (\delta (q',\sigma _{q'})\in
  F)^{1-\epsilon _{q'}}\) holds (where as before \(\phi ^0\) denotes \(\neg\phi
  \) and \(\phi ^1\) denotes \(\phi \)). If \(\theta \in Q\) then \(\theta
  = q\) if and only if
  \[\bigwedge _{q'\ne q} (\delta (\theta ,\sigma _{q'})\in F)^{\epsilon _{q'}}
  \]
  holds.
  Consider then the following formula in the variables \(\overline{x} = (x_{q'}
  :  q'\ne q)\) and \(y\):
  \[\phi (\overline{x};y) = \bigwedge _{q'\ne q} (x_{q'} + y \in A)^{\epsilon
    _{q'}}
  \]
  We show that \(\phi \) is unstable in \(\mathbb{N}\), and hence since \(\phi
  \) is a Boolean combination of instances of \(x+y\in A\) that \(A\) is
  unstable in \(\mathbb{N}\).

  Recall that \(L_q\) satisfies the hypotheses of
  \cref{lemma:nongenspecialcase}; so for some \(K<\omega \) there is an
  \(N\)-ladder \((d_i,e_i : i \le N)\) for \(x+_Ky\in L_q\).  We may assume
  each \(\abs{d_i} = K\). Take any \(\mu \in\Sigma ^*\) such that \(\delta
  (q_0,\mu ) = q\), and let
  \begin{ea}
    b_{i,q'} &=& [\mu d_i \sigma _{q'}] \\
    c_i &=& [0^{\abs\mu }e_i]
  \end{ea}
  These \(\overline{b_i} := (b_{i,q'}: q'\ne q),c_i\) will be our ladder for
  \(\phi \).  Note that \(b_{i,q'} + c_j = [\mu (d_i +_K e_j) \sigma _{q'}]\).
  Then
  \begin{ea}
    \phi (\overline{b_i}; c_j) 
    &\iff& \bigwedge _{q'\ne q}(b_{i,q'} + c_j\in A)^{\epsilon _{q'}} \\
    &\iff& \bigwedge _{q'\ne q}(\delta (q_0,\mu  (d_i +_K e_j) \sigma
    _{q'})\in F)^{\epsilon _{q'}} \\
    &\iff& \bigwedge _{q'\ne q}(\delta (\delta (q,d_i +_K e_j), \sigma
    _{q'})\in F)^{\epsilon _{q'}} \\
    &\iff& \delta (q,d_i+_Ke_j)= q \\
    &\iff& d_i+_Ke_j\in L_q \\
    &\iff& i\le j.
  \end{ea}
  So \(\phi \) is unstable in \(\mathbb{N}\), and thus \(A\) is unstable in
  \(\mathbb{N}\).
\end{lrpf}
We can now do the case \(A\subseteq \mathbb{Z}\):
\begin{lrpf}{thm:stablenongen}
  Suppose \(A\subseteq \mathbb{Z}\) is \(d\)-automatic but neither \(d\)-sparse
  nor generic in \(\mathbb{Z}\).
  \begin{caselist}
    \case
      Suppose one of \(A\cap\mathbb{N}\) and \(-A\cap\mathbb{N}\) is generic in
      \(\mathbb{N}\) and the other is \(d\)-sparse. Then taking finitely many
      translates and unioning we get a set \(B\) where (say)
      \(B\cap\mathbb{N}\) is \(d\)-sparse and \(B\supseteq -\mathbb{N}\). (Note
      that \(d\)-sparsity is closed under translation and finite union.) Recall
      by \cite[Theorem 3.8]{yu97} that the set of prefixes of a sparse set is
      also sparse, and in particular that every sparse set has a forbidden
      prefix.  So there is \(\sigma \in\Sigma ^*\) such that \(\sigma \) is not
      a prefix of any canonical representative of an element of
      \(B\cap\mathbb{N}\); by possibly appending a \(1\), we may assume that
      \(\sigma \) has no trailing zeroes. So if \(r=[\sigma ]\) and
      \(s=d^{\abs\sigma }\) then \((r + s\mathbb{N})\cap B = \emptyset \); so
      \((r + s\mathbb{Z})\cap B = r+s\mathbb{Z}_{<0}\), and thus
      \((r+s\mathbb{Z})\cap B\) is unstable in \((\mathbb{Z},+)\) since if
      \(x,y\in s\mathbb{Z}\) then
      \[x\le y\iff r+x-y -s\in r+s\mathbb{Z}_{<0}
      \]
      So \((x+y\in r+s\mathbb{Z})\wedge (x+y\in B)\) is unstable in
      \((\mathbb{Z},+)\). But \(x+y\in r+s\mathbb{Z}\) is stable in
      \((\mathbb{Z},+)\), since it's definable in \((\mathbb{Z},+)\); so \(B\)
      is unstable in \((\mathbb{Z},+)\). So since \(B\) is a finite union of
      translates of \(A\) we get that \(A\) is unstable in \((\mathbb{Z},+)\).
    \case
      Suppose otherwise.  Since \(A\) isn't generic in \(\mathbb{Z}\), at most
      one of \(A\cap\mathbb{N}\) or \(-A\cap\mathbb{N}\) is generic in
      \(\mathbb{N}\); likewise with \(d\)-sparse.  Since we precluded the
      previous case we know there can't be one of each, and generic in
      \(\mathbb{N}\) and \(d\)-sparse are contradictory. So one of
      \(A\cap\mathbb{N}\) or \(-A\cap\mathbb{N}\) is neither generic in
      \(\mathbb{N}\) nor \(d\)-sparse. Note that \(A\) is stable in
      \((\mathbb{Z},+)\) if and only if \(-A\) is. Hence replacing \(A\) by
      \(-A\) if necessary we may assume \(A\cap \mathbb{N}\) is neither generic
      in \(\mathbb{N}\) nor \(d\)-sparse.  Then by \cref{prop:stablenatnongen}
      there are arbitrarily large ladders in \(\mathbb{N}\) for \(x+y\in
      A\cap\mathbb{N}\); since \(\mathbb{N}\) is closed under addition, we get
      that these are also ladders in \(\mathbb{Z}\) for \(x+y\in A\). Hence
      \(A\) is unstable in \((\mathbb{Z},+)\).
      \qedhere
  \end{caselist}
\end{lrpf}
As an illustration of our theorem we note that the following automatic sets are
not stable in \((\mathbb{Z},+)\). Indeed, it is easily checked that they are
all neither sparse nor generic.
\begin{corollary}
  The following automatic sets are unstable in \((\mathbb{Z},+)\):
  \begin{itemize}
    \item
      The set of \(a\in\mathbb{Z}\) such that the canonical base-\(d\)
      representation of \(a\) ends in \(\pm1\) (assuming \(d>2\)).
    \item
      The set of \(a\in\mathbb{Z}\) such that the canonical base-\(d\)
      representation of \(a\) doesn't contain a \(0\) (assuming \(d>2\)).
    \item
      The set of \(a\in\mathbb{Z}\) such that the canonical base-\(d\)
      representation of \(a\) is of even length.
    \item
      The set of \(a\in\mathbb{Z}\) such that in the canonical binary
      representation of \(a\) takes the form \(0^{k_0}10^{k_1}1\cdots
      10^{k_m}1\) or \(0^{k_0}(-1)0^{k_1}(-1)\cdots (-1)0^{k_m}(-1)\) for some
      even \(k_0,\ldots ,k_m\) (possibly zero); i.e.\ does not contain a block
      of zeroes of odd length.  These are precisely the \(a\in \mathbb{Z}\)
      such that the Baum-Sweet sequence has a \(1\) in the \(\abs a\th\)
      position. See \cite[Section 5.1]{allouche03} for more details on the
      Baum-Sweet sequence.
  \end{itemize}
\end{corollary}

The converse of \cref{thm:stablenongen} is certainly false. For example, let
\(A\subseteq \mathbb{Z}\) be as in the example at the beginning of
\cref{sec:sparsestable}; so \(A\) is \(d\)-sparse and unstable in
\((\mathbb{Z},+)\). Then the complement of \(A\) remains unstable, and is
generic since \(A\) doesn't contain a pair of adjacent integers.

\section{The general case}
Gabriel Conant pointed out to me in private communications that
\cref{thm:charstablesparse,thm:stablenongen}, together with \cite[Theorem 2.3
(iv)]{conant20pt}, allow us deal with arbitrary \(d\)-automatic stable subsets
of \(\mathbb{Z}\).
\begin{theorem}
  \label{thm:charstableauto}
  Suppose \(A\subseteq \mathbb{Z}\) is \(d\)-automatic and stable in
  \((\mathbb{Z},+)\). Then \(A\) is a finite Boolean combination of
  \begin{itemize}
    \item
      cosets of subgroups of \((\mathbb{Z},+)\), and
    \item
      translates of finite sums of sets of the form \(C(a;\delta )\).
  \end{itemize}
\end{theorem}
\begin{rpf}
  It is known that stable subsets of a group are close to being a finite union
  of cosets, in the sense that they have non-generic symmetric difference with
  such; see \cite[Theorem 2.3 (iv)]{conant20pt} (taking \(\delta (x,y)\) to be
  \(x+y\in A\) and \(\phi (x)\in\operatorname{Def}_\delta (G)\) to be \(x\in
  A\)). So there is a subgroup \(H\le \mathbb{Z}\) and a union \(Y\) of cosets
  of \(H\) such that \(Z := A\symdif Y\) is non-generic in \(\mathbb{Z}\).
  Since \(Y\) is a union of cosets it is also stable in \((\mathbb{Z},+)\) and
  \(d\)-automatic. Hence \(Z\) is both \(d\)-automatic and stable in
  \((\mathbb{Z},+)\). \Cref{thm:stablenongen} yields that \(Z\) is
  \(d\)-sparse, and then \cref{thm:charstablesparse} yields that \(Z\) is a
  finite Boolean combination of translates of finite sums of \(C(a;\delta)\).
  Hence \(A = Z\symdif Y\) is a finite Boolean combination of sets of the
  desired form.
\end{rpf}
\begin{corollary}
  \label{cor:autosparsechar}
  Suppose \(A\subseteq \mathbb{Z}\) is \(d\)-automatic.  \Tfae{}:
  \begin{enumerate}
    \item
      \(\thy(\mathbb{Z},+,A)\) is stable.
    \item
      \(A\) is stable in \((\mathbb{Z},+)\).
    \item
      \(A\) is a finite Boolean combination of
      \begin{itemize}
        \item
          cosets of subgroups of \((\mathbb{Z},+)\), and
        \item
          translates of finite sums of sets of the form \(C(a;\delta )\).
      \end{itemize}
    \item
      \(A\) is definable in \((\mathbb{Z},\mathcal{F})\).
  \end{enumerate}
\end{corollary}
\begin{rpf}
  That (1) \(\implies \) (2) and (3) \(\implies \) (4) is clear. That (2)
  \(\implies \) (3) is \cref{thm:charstableauto}. Finally, that (4) \(\implies
  \) (1) is \cite[Theorem A]{moosa04}.
\end{rpf}

\section{Two NIP expansions of \texorpdfstring{\((\mathbb{Z},+)\)}{(Z,+)}}
In this final section we show how to apply automata-theoretic methods to
produce some NIP expansions of \((\mathbb{Z},+)\); see \cite{simon15} for
background on NIP.  

\subsection{\texorpdfstring{\((\mathbb{Z},+,<,d^\mathbb{N})\)}
  {(Z,+,<,d\textasciicircum N)} is NIP}

Fix \(d>0\). That \(\thy(\mathbb{Z},+,<,d^\mathbb{N})\) is NIP was shown
recently by Lambotte and Point (it is an instance of \cite[Corollary
2.33]{lambotte20}), but our proof is novel and short.  It will be convenient to
work in \((\mathbb{N},+)\) rather than \((\mathbb{Z},+,<)\).  Since
\((\mathbb{Z},+,<,d^\mathbb{N})\) is interpretible in
\((\mathbb{N},+,d^\mathbb{N})\), it will suffice to prove:
\begin{theorem}
  \label{thm:presbpowersnip}
  \(\thy(\mathbb{N},+, d^\mathbb{N})\) is NIP.
\end{theorem}
Before proving the theorem, let us observe that since all \(d\)-sparse subsets
of \(\mathbb{N}\) are definable in \((\mathbb{N},+,d^\mathbb{N})\)---see
\cite[Theorem 5]{semenov79}---and as \(A\subseteq \mathbb{Z}\) is \(d\)-sparse
if and only if both \(A\cap \mathbb{N}\) and \(-A\cap\mathbb{N}\) are, we get:
\begin{corollary}
  The expansion of \((\mathbb{Z},+)\) by all \(d\)-sparse subsets is NIP.
\end{corollary}
Our proof of \cref{thm:presbpowersnip} will make use of a result of Chernikov
and Simon on NIP pairs of structures; we briefly recall their setup and result.
We let \(L = \set{+}\) and \(\mathcal{N} = (\mathbb{N},+)\); we fix
\(\thy(\mathcal{N})\) as our ambient theory.
\begin{definition}
  Let \(L_P\) be \(L\) expanded by a unary predicate \(P\). A \emph{bounded
  \(L_P\)-formula} is one of the form \((Q_1x_1\in P)\cdots (Q_nx_n\in P) \phi
  \) for some quantifiers \(Q_i\) and some \(\phi \in L\). If \(M\) is an
  \(L\)-structure and \(A\subseteq M\) we say \(A\) is \emph{bounded in \(M\)}
  if every \(L_P\)-formula is \(\thy(M,A)\)-equivalent to a bounded one.
\end{definition}
\begin{definition}
  Suppose \(M\) is a structure and \(A\subseteq M\). The \emph{induced
  structure} \(A_M\) of \(M\) on \(A\) has domain \(A\) and atomic relations
  \(D\cap A^n\) for each \(\emptyset \)-definable \(D\subseteq M^n\).
\end{definition}
\begin{fact}[\cite[Corollary 2.5]{chernikov13}]
  \label{fact:chernikovsimon}
  Suppose \(M\) is a structure and \(A\subseteq M\) is bounded in \(M\). If
  \(\thy(M)\) and \(\thy(A_M)\) are NIP then so is \(\thy(M,A)\).
\end{fact}
We wish to apply this to \((M,A) = (\mathcal{N},d^\mathbb{N})\). Boundedness
follows from earlier work of Point:
\begin{proposition}
  \label{prop:powersbounded}
  \(d^\mathbb{N}\) is bounded in \(\mathcal{N}\).
\end{proposition}
\begin{rpf}
  \cite[Propositions 9 and 11]{point00} say that
  \(\thy(\mathbb{N},+,\dot-,<,0,1,\frac\cdot n,\lambda _d,S,S^{-1})_{n\ge1}\)
  admits quantifier elimination, where
  \begin{itemize}
    \item
      \(S(d^n) = d^{n+1}\) and \(S(a) = a\) for other \(a\)
    \item
      \(S^{-1}(d^{n+1}) = d^n\) and \(S^{-1}(a) = a\) for other \(a\)
    \item
      \(\lambda _d(x) = d^{\floor{\log_d(x)}}\) for \(x>0\) and \(\lambda _d(0)
      = 0\).
  \end{itemize}
  It then remains to show that any quantifier-free formula in this signature is
  equivalent to a bounded \(L_P\)-formula.  But a quantifier-free formula
  \(\phi (\ldots ,\lambda _d(t), \ldots )\) involving \(\lambda _d\) is
  equivalent to
  \[(\exists x\in d^\mathbb{N})\pars[\Big]{(x\le t)\wedge (\forall y\in
      d^\mathbb{N})\neg(x <
    y \le t)\wedge \phi (\ldots ,x,\ldots )}\vee \pars[\Big]{(t=0)\wedge \phi
    (\ldots ,0,\ldots )}
  \]
  So at the cost of quantifying over \(d^\mathbb{N}\) we can eliminate
  occurrences of \(\lambda _d\); we can similarly dispense with occurrences of
  \(S\) and \(S^{-1}\). Repeatedly applying this yields that any
  quantifier-free formula in the given signature is equivalent to one of the
  form \((Q_1x_1\in d^\mathbb{N})\cdots (Q_nx_n\in d^\mathbb{N})\psi \) where
  \(\psi \) is a formula in \(\set{0,1,+,\dot-,\frac\cdot n}_{n\ge1}\).  But
  since \((\mathbb{N},0,1,+,\dot-,\frac\cdot n)_{n\ge1}\) is a definitional
  expansion of \((\mathbb{N},+)\), we get that \(\phi \) is equivalent to a
  bounded \(L_P\)-formula.
\end{rpf}
It is well-known that \(\mathcal{N}\) is NIP; it is definable in
\((\mathbb{Z},+,<)\), which is NIP as all ordered abelian groups are (see
\cite{gurevich84}).  It remains to show that the induced structure
\((d^\mathbb{N})_{\mathcal{N}}\) is NIP. 

The following is well-known; see e.g.\ \cite[Theorem 6.1]{bruyere94}, of which
it is a weakening.
\begin{fact}
  \label{fact:dfblauto}
  All definable subsets of \(\mathcal{N}\) are \(d\)-automatic.
\end{fact}
We therefore wish for a description of how \(d\)-automatic sets can intersect
\(d^\mathbb{N}\).
\begin{proposition}
  \label{prop:powerspresbdfbl}
  If \(X\subseteq \mathbb{N}^n\) is \(d\)-automatic then the relation
  \[\set{(k_1,\ldots ,k_n) \in\mathbb{N}^n: (d^{k_1},\ldots ,d^{k_n})\in X}
  \]
  is definable in \((\mathbb{N},+)\).
\end{proposition}
\begin{rpf}
  % By symmetry and disjunction it suffices to give a first-order definition in
  % the case \(k_1\le \cdots \le k_n\). We in fact show that given a
  % deterministic finite automaton \((\set{0,\ldots ,d-1} ^n,Q, q_0,\delta , F)\)
  % and any \(q_1,q_2\in Q\) that
  % \[\set*{(k_1,\ldots ,k_n)\in\mathbb{N}^n : \delta ^*\pars*{q_1, \fcol{0^{k_1}1
  %   0^{k_n-k_1}}\vdots {0^{k_{n-1}}10^{k_n-k_{n-1}}}{ 0^{k_n}1}} = q_2}
  % \]
  % is definable in \((\mathbb{N},+)\). (Here \(0^k\) denotes the string of \(k\)
  % zeroes.) We apply induction on \(n\); the case \(n=0\) is immediate.
  % For the induction step, note that since there are finitely many states,
  % eventually
  % \[\delta ^*\pars*{q_1, \tcol{0^\ell }\vdots {0^\ell }}
  % \]
  % is cyclic in \(\ell \).
  By symmetry and disjunction it suffices to check the case
  \(k_1\le\cdots \le k_n\).

  It will be more convenient to work with \((X\cap\mathbb{N}_{>0}^n)-1\), which
  is also \(d\)-automatic.  Then taking
  \[\sigma _i = \begin{pmatrix}0\\\vdots
    \\0\\d-1\\\vdots \\d-1\end{pmatrix}
  \] 
  with \(i-1\) zeroes, we get for \(k_1\le\cdots \le k_n\) that
  \[(d^{k_1},\ldots ,d^{k_n})\in X
    \iff (d^{k_1}-1,\ldots ,d^{k_n}-1)\in (X\cap\mathbb{N}_{>0}^n)-1
    \iff [\sigma _1^{k_1}\sigma _2^{k_2-k_1}\cdots \sigma _n^{k_n-k_{n-1}}] \in
    (X\cap\mathbb{N}_{>0}^n)-1
  \]
  (since the base-\(d\) representation of \(d^{k_i}-1\) consists of \(d-1\)
  repeated \(k_i\) times). But by \cref{prop:sparseintersregpresb} the last
  condition is definable in \((\mathbb{N},+)\), as desired.
\end{rpf}
Our theorem now follows easily:
\begin{lrpf}{thm:presbpowersnip}
  \Cref{prop:powerspresbdfbl,fact:dfblauto} imply that the map \(k\mapsto d^k\)
  induces an interpretation of \((d^\mathbb{N})_\mathcal{N}\) in
  \((\mathbb{N},+)\). But \(\thy(\mathbb{N},+)\) is NIP; so
  \(\thy(d^\mathbb{N})_\mathcal{N}\) is NIP. But \(d^\mathbb{N}\) is bounded in
  \(\mathcal{N}\) by \cref{prop:powersbounded}, and \(\mathcal{N}\) is NIP. So
  \(\thy(\mathcal{N},d^\mathbb{N})\) is NIP by \cref{fact:chernikovsimon}.
\end{lrpf}

\subsection{\texorpdfstring{\((\mathbb{Z},+,d^\mathbb{N},\dpowmult)\)}
  {(Z,+,d\textasciicircum N,x|d\textasciicircum N)} is NIP}

Next we consider the expansion of \((\mathbb{Z},+)\) by the monoid
\((d^\mathbb{N},\times )\). Note that as the ordering on \(d^\mathbb{N}\) is
definable here, \((\mathbb{Z},+,d^\mathbb{N},\dpowmult)\) is not stable.
However:
\begin{theorem}
  \label{thm:multpowersnip}
  \(\thy(\mathbb{Z},+,d^\mathbb{N},\dpowmult)\) is NIP.
\end{theorem}
Surprisingly, our methods apply even though \(\dpowmult\) itself isn't
\(d\)-automatic: since
\[\bracks*{\tcol{0^i1}{0^i1}{0^{i+1}}\cdot
  \tcol{0^{j+1}}{0^{j+1}}{0^j1}}\in \dpowmult
  \iff
  i = j+1
\]
it follows from the Myhill-Nerode theorem (see e.g.\ \cite[Theorem 4.7]{yu97})
that the set of canonical representations of elements of \(\dpowmult\) isn't
regular.  The reason automatic methods still apply is \cref{fact:dfblauto},
together with the following generalization of \cref{prop:powerspresbdfbl},
which tells us that the interaction between iterated concatenation and
membership in automatic sets can be described using Presburger arithmetic.
\begin{lemma}
  \label{lemma:autopresdfbl}
  Suppose \(X\subseteq \mathbb{Z}^m\) is \(d\)-automatic and \((\ell
  _{11},\ldots ,\ell _{1n_1}),\ldots ,(\ell _{m1},\ldots ,\ell _{mn_m})\) are
  tuples from \(\Sigma _\pm\). Then the relation
  \[\set*{(k_{ij}) : 
      \tcol{\relax[\ell _{11}^{k_{11}}\cdots \ell _{1n_1}^{k_{1n_1}}]}
      {\vdots }
      {\relax[\ell _{m1}^{k_{m1}}\cdots \ell _{mn_m}^{k_{mn_m}}]}
      \in X
    }
    \subseteq \mathbb{N}^{n_1}\times \cdots \mathbb{N}^{n_m}
  \]
  is definable in \((\mathbb{N},+)\).
\end{lemma}
\begin{rpf}
  We show that for any \(\ell _{ij}\), any automaton \((Q,q_0,\delta ,F)\),
  and any \(q_1,q_2\in Q\) the relation
  \[\set*{(k_{ij}) : \delta \pars*{q_1,P\tcol{
        \ell _{11}^{k_{11}}\cdots \ell _{1n_1}^{k_{1n_1}}
      }{\vdots }{
        \ell _{m1}^{k_{m1}}\cdots \ell _{mn_m}^{k_{mn_m}}
      }
    } = q_2}
  \]
  is definable in \((\mathbb{N},+)\), where \(P\colon (\Sigma _\pm^*)^m\to
  (\Sigma _\pm^m)^*\) takes in \(m\) strings and pads them on
  the right with zeroes so they all have the same length as the longest one.
  This claim, applied to an automaton for the set of representations over
  \(\Sigma _\pm^m\) of elements of \(X\), yields the desired result.

  We apply induction on \((m,n_1,\ldots ,n_m)\).  The base case \(m=0\) is
  vacuous.  For the induction step, suppose first that some \(n_i = 0\); say
  for ease of notation that \(i=1\). Then we can construct an automaton
  \((Q,q_0,\delta ',F)\) over \(\Sigma _\pm^{m-1}\) such that
  \(\delta' (q,\sigma ) = \delta \pars*{q,\col{0^{\abs\sigma }}\sigma }\) for
  any \(\sigma \in\Sigma _\pm^{m-1}\); that is, it behaves like the
  original automaton would if the input had an extra string of zeroes attached.
  In particular we have
  \[\delta \pars*{q_1, P\fcol{\epsilon }{
        \ell _{21}^{k_{21}}\cdots \ell _{2n_2}^{k_{2n_2}}
      }\vdots {
        \ell _{m1}^{k_{m1}}\cdots \ell _{mn_m}^{k_{mn_m}}
    }} = q_2
    \iff
    \delta '\pars*{q_1, P\tcol{
        \ell _{21}^{k_{21}}\cdots \ell _{2n_2}^{k_{2n_2}}
      }\vdots {
        \ell _{m1}^{k_{m1}}\cdots \ell _{mn_m}^{k_{mn_m}}
    }} = q_2
  \]
  and by the induction hypothesis the latter is definable in
  \((\mathbb{N},+)\).

  Suppose then that no \(n_i=0\). Suppose \(k_{11}\) is minimum among the
  \(k_{i1}\). Then if
  \[q = \delta \pars*{q_1, \tcol{\ell _{11}^{k_{11}}}{\vdots }{\ell
    _{m1}^{k_{11}}}}
  \]
  then our relation is equivalent to
  \[\delta \pars*{q, P\fcol{
      \hphantom{\ell _{11}^{k_{11}-k_{11}}}\ell _{12}^{k_{12}}\cdots \ell
      _{1n_1}^{k_{1n_1}}
    }{
      \ell _{21}^{k_{21}-k_{11}}\ell _{22}^{k_{22}}\cdots
      \ell _{2n_2}^{k_{2n_2}}
    }\vdots {
      \ell _{m1}^{k_{m1}-k_{11}}\ell _{m2}^{k_{m2}}\cdots
      \ell _{mn_m}^{k_{mn_m}}
    }} = q_2
  \]
  which by the induction hypothesis is definable in \((\mathbb{N},+)\).
  Similarly we get definability in the case \(k_{i1}\) is minimum for some
  \(i>1\). So taking disjunctions we get that our relation is definable in
  \((\mathbb{N},+)\).
\end{rpf}
For our proof of \cref{thm:multpowersnip} it will be convenient to assume
\(d\ge8\). In fact this suffices: consider for example the case \(d=4\).
Assuming the theorem holds when \(d=4^2=16\), we get that
\((\mathbb{Z},+,16^\mathbb{N},\powmult{16})\) is NIP.  But \(\powmult{4}\) is
definable in \((\mathbb{Z},+,16^\mathbb{N},\powmult{16})\): we have
\((a,b,c)\in \powmult{4}\) if and only if \((4^ia,4^jb, 4^{i+j}c)\in
\powmult{16}\) for some \(i,j\in\set{0,1}\). This is because \(x\) is a power
of \(4\) if and only if one of \(x,4x\) is a power of \(16\). So
\((\mathbb{Z},+,4^\mathbb{N},\powmult{4})\) is definable in
\((\mathbb{Z},+,16^\mathbb{N},\powmult{16})\), and is thus NIP. Similar
arguments work for all \(2\le d< 8\).
\begin{lrpf}{thm:multpowersnip}
  We assume \(d\ge8\). We will apply an extension due to Conant and Laskowski
  of the result of Chernikov and Simon we used previously
  (\cref{fact:chernikovsimon}). Since these results only apply to subsets of
  the domain, our first task is to encode \(d^\mathbb{N}\) and \(\dpowmult\) as
  such.  Let
  \[B = d^\mathbb{N}\cup\set{[7^i6^j4^i] : i,j\in\mathbb{N}}
  \]
  The point is that from \(B\) we will be able to extract both \(d^\mathbb{N}\)
  and
  \[\set*{\frac{a-1}{d-1} + 2\frac{b-1}{d-1} + 4\frac{c-1}{d-1} : (a,b,c)\in
    \dpowmult,a \le b}
  \]
  These together will be enough to recover \(\dpowmult\).
  \begin{claim}
    \(d^\mathbb{N}\) and \(\dpowmult\) are definable in \((\mathbb{Z},+,B)\).
  \end{claim}
  \begin{rpf}
    Note first that \(d^\mathbb{N}\) is definable in \((\mathbb{Z},+,B)\): we
    have \(a\in d^\mathbb{N}\) if and only if \(a = 1\) or \(0\ne a\in B\) and
    \(a\equiv 0\pmod d\). I now claim that \((a,b,c)\in \dpowmult\) with \(a\le
    b\) if and only if \(a,b,c\in d^\mathbb{N}\) and
    \[\frac{a-1}{d-1} + 2\frac{b-1}{d-1} + 4\frac{c-1}{d-1} \in B
    \]
    For the left-to-right direction, note that if \((d^i,d^j,d^{i+j})\in
    \dpowmult\)
    with \(i\le j\) then
    \[\frac{d^i-1}{d-1} + 2\frac{d^j-1}{d-1} + 4\frac{d^{i+j}-1}{d-1}
      = [1^i] + [2^j] + [4^{i+j}]
      = [7^i6^{j-i}4^i] \in B
    \]
    For the right-to-left direction, suppose \(d^i,d^j,d^k\) satisfy
    \[
      [1^i] + [2^j] + [4^k] =
      \frac{d^i-1}{d-1} + 2\frac{d^j-1}{d-1} + 4\frac{d^k-1}{d-1}
      \in B
    \]
    If \(i=j=k=0\) then \((d^i,d^j,d^k)\in \dpowmult\) and \(d^i\le d^j\), as
    desired; suppose then that at least one is non-zero. Then \([1^i]+[2^j] +
    [4^k]\not\equiv 0\pmod d\), so \([1^i]+[2^j] + [4^k]\in B\setminus
    d^\mathbb{N}\), and is thus equal to \([7^{i'}6^{j'}4^{i'}] = [1^{i'}] +
    [2^{j'+i'}] + [4^{2i'+j'}]\) for some \(i',j'\).

    But the map \((x,y,z)\mapsto [1^x] + [2^y] + [4^z]\) is
    injective. Indeed, we can represent \([1^x] + [2^y] + [4^z]\) by element of
    \(\set{1,\ldots ,7}^*\); note that each element of \(\set{1,\ldots ,7}\)
    can be represented uniquely as a sum of a subset of \(\set{1,2,4}\). We can
    then recover \(x\) from the canonical representation of
    \([1^x]+[2^y]+[4^z]\) as the number of occurrences of \(\ell  \in
    \set{1,\ldots 7}\) that use a \(1\) in this sum representation; we can
    likewise recover \(y,z\).

    So since \([1^i] + [2^j] + [4^k] = [1^{i'}] + [2^{j'+i'}] + [4^{2i'+j'}]\)
    we get by injectivity that \(j = j'+i'\ge i' = i\) and \(k = 2i'+j' = i +
    j\); so \((d^i,d^j,d^k)\in \dpowmult\) and \(d^j\ge d^i\), as desired.

    But \((a,b,c)\in \dpowmult\iff (b,a,c)\in \dpowmult\); so
    \[(x\le y\wedge (x,y,z)\in \dpowmult)\vee (y\le x\wedge  (y,x,z)\in
      \dpowmult)
    \]
    defines \(\dpowmult\) in \((\mathbb{Z},+,B)\).
  \end{rpf}
  So it suffices to show that \((\mathbb{Z},+,B)\) is NIP. We again check that
  the induced structure on \(B\) is NIP. When using \cref{fact:chernikovsimon},
  we only concerned ourselves with the structure induced from the \(\emptyset
  \)-definable sets; however, to use the result of Conant and Laskowski, we
  will need that the structure induced by all sets definable with parameters
  from \(\mathbb{Z}\) is NIP.
  \begin{claim}
    \label{claim:inducednip}
    Let \(\mathcal{Z}\) be \((\mathbb{Z},+)\) expanded by names for all the
    constants. Then the induced structure \(B_\mathcal{Z}\)~is~NIP.
  \end{claim}
  \begin{rpf}
    Let \(D = \set{(e_1,1,0,0) : e_1\in\mathbb{N}} \cup \set{(0,0,e_3,e_4):
    e_3,e_4\in\mathbb{N} }\subseteq \mathbb{N}^4\); note that \(D\) is
    definable in \((\mathbb{N},+)\). Consider \(\Phi \colon \mathbb{N}^4\to
    \mathbb{Z}\) given by \((e_1,e_2,e_3,e_4)\mapsto
    [0^{e_1}1^{e_2}7^{e_3}6^{e_4}4^{e_3}]\); note that \(\Phi (D)\subseteq B\),
    and in fact \(\Phi \colon D\to B\) is bijective.  I claim that \(\Phi \)
    defines an interpretation of \(B_\mathcal{Z}\) in \((\mathbb{N},+)\).
    Recall that \((\mathbb{Z},+,0,1,\delta \mathbb{N})_{ \delta >0}\) admits
    quantifier elimination (see e.g.\ \cite[Exercise
    3.4.6]{marker06}).  So if \(X\subseteq \mathbb{Z}\) is
    definable in \(\mathcal{Z}\) then \(X\) is a Boolean combination of
    congruences and equalities, and hence \(X\cap\mathbb{N}\) is definable in
    \((\mathbb{N},+)\); likewise with \(-X\cap\mathbb{N}\). So since
    \(\mathbb{N}\) is a \(d\)-automatic subset of \(\mathbb{Z}\) and
    \(d\)-automatic sets are closed under Boolean combinations we get that
    \(X\) is \(d\)-automatic. One argues similarly that if \(X\subseteq
    \mathbb{Z}^m\) is definable in \(\mathcal{Z}\) then \(X\) is
    \(d\)-automatic. So to show that \(\Phi \) defines an interpretation it
    suffices to show that whenever \(X\subseteq \mathbb{Z}^m\) is
    \(d\)-automatic we have that
    \[\set*{(e_{ij}) \in D^m:
        \tcol{\relax[0^{e_{11}}1^{e_{12}}7^{e_{13}}6^{e_{14}}4^{e_{13}}]}\vdots
      {\relax[0^{e_{m1}}1^{e_{m2}}7^{e_{m3}}6^{e_{m4}}4^{e_{m3}}]} \in X}
    \]
    is definable in \((\mathbb{N},+)\). But this follows from
    \cref{lemma:autopresdfbl} (and definability of \(D\)). So \(\Phi \) defines
    an interpretation of \(B_\mathcal{Z}\) in \((\mathbb{N},+)\); so
    \(B_\mathcal{Z}\) is NIP.
  \end{rpf}
  Now by \cite[Theorem 2.9]{conant20} we get since \(\thy(\mathbb{Z},+)\) is
  weakly minimal (see e.g.\ \cite[Proposition 3.1]{conant20}) and
  \(B_\mathcal{Z}\) is NIP that \((\mathbb{Z},+,B)\) is NIP.  So
  \((\mathbb{Z},+,d^\mathbb{N},\dpowmult)\) is NIP.
\end{lrpf}
Despite the similarity of methods in
\cref{thm:presbpowersnip,thm:multpowersnip}, we don't know whether
\(\thy(\mathbb{Z},+,<,d^\mathbb{N},\dpowmult)\) is NIP. One might hope to apply
\cref{fact:chernikovsimon} with \((\mathbb{Z},+,<)\) as the base NIP structure
and \(B\) as the new predicate.  Indeed, as in the proof of
\cref{claim:inducednip} one can show that the induced structure on \(B\) is NIP
by observing that the definable subsets of \((\mathbb{Z},+,<)\) are
\(d\)-automatic.  Checking boundedness, however, isn't simply a matter of
adapting the arguments of \cref{thm:presbpowersnip} as the quantifier
elimination result of Point that applied to \(d^\mathbb{N}\) doesn't seem to
apply to \(B\). Nor does the result of Conant and Laskowski yield boundedness
as \((\mathbb{Z},+,<)\) is not weakly minimal.
% we used a quantifier elimination result of Point in \cite{point00} to check
% boundedness: that if \(R\subseteq \mathbb{N}\) is enumerated by an
% \emph{almost sparse} sequence then an appropriate expansion of
% \((\mathbb{N},+,R)\) admits quantifier elimination. Unfortunately, the \(B\)
% we used above isn't almost sparse: it isn't ultimately periodic modulo \(d\).
% \begin{todo} Is that what axiom (5) says? Sure doesn't look like it. Take \(E
% = 0\).  \end{todo}
So if one wishes to use our approach to show that
\(\thy(\mathbb{Z},+,<,d^\mathbb{N},\dpowmult)\) is NIP one needs a new way to
check boundedness.

One can restate \cref{thm:multpowersnip} as saying that expanding
\((\mathbb{Z},+)\) by a singly generated submonoid of \((\mathbb{Z}\setminus
\set0,\times )\) yields an NIP structure. It would be natural to ask about
finitely generated submonoids in general, but it
seems unlikely that our automata-theoretic methods will apply as there is no
obvious choice of \(d\) in general.
\printbibliography
\end{document}